\input amstex
\catcode`\"=12

\font\black=cmbx10
\font\sblack=cmbx7
\font\ssblack=cmbx5
\font\blackital=cmmib10  \skewchar\blackital='177
\font\sblackital=cmmib7  \skewchar\sblackital='177
\font\ssblackital=cmmib5  \skewchar\ssblackital='177
\font\sanss=cmss10

\font\teneurm=eurm10

\font\blackboard=msbm10
\font\sblackboard=msbm7
\font\ssblackboard=msbm5

\font\fraktur=eufm10

\newfam\blfam
\textfont\blfam=\black
\scriptfont\blfam=\sblack
\scriptscriptfont\blfam=\ssblack

\newfam\bifam
\textfont\bifam=\blackital
\scriptfont\bifam=\sblackital
\scriptscriptfont\bifam=\ssblackital

\newfam\bbfam
\textfont\bbfam=\blackboard
\scriptfont\bbfam=\sblackboard
\scriptscriptfont\bbfam=\ssblackboard
\def\bb#1{{\fam\bbfam\relax#1}}

\newfam\gpfam
\def\ss#1{\hbox{$\fam\gpfam\relax#1\textfont\gpfam=\sanss$}}

\def\eu#1{\hbox{$\fam\gpfam\relax#1\textfont\gpfam=\teneurm$}}

\def\fr#1{\hbox{$\fam\gpfam\relax#1\textfont\gpfam=\fraktur$}}

\mathchardef\za="710B  
\mathchardef\zb="710C  
\mathchardef\zg="710D  
\mathchardef\zd="710E  
\mathchardef\zve="710F 
\mathchardef\zz="7110  
\mathchardef\zh="7111  
\mathchardef\zvy="7112 
\mathchardef\zi="7113  
\mathchardef\zk="7114  
\mathchardef\zl="7115  
\mathchardef\zm="7116  
\mathchardef\zn="7117  
\mathchardef\zx="7118  
\mathchardef\zp="7119  
\mathchardef\zr="711A  
\mathchardef\zs="711B  
\mathchardef\zt="711C  
\mathchardef\zu="711D  
\mathchardef\zvf="711E 
\mathchardef\zq="711F  
\mathchardef\zc="7120  
\mathchardef\zw="7121  
\mathchardef\ze="7122  
\mathchardef\zy="7123  
\mathchardef\zf="7124  
\mathchardef\zvr="7125 
\mathchardef\zvs="7126 
\mathchardef\zf="7127  
\mathchardef\zG="7000  
\mathchardef\zD="7001  
\mathchardef\zY="7002  
\mathchardef\zL="7003  
\mathchardef\zX="7004  
\mathchardef\zP="7005  
\mathchardef\zS="7006  
\mathchardef\zU="7007  
\mathchardef\zF="7008  
\mathchardef\zW="700A  

\catcode`\"=\active
\TagsOnRight
\documentstyle{amsppt}
\hsize=33pc
\hoffset=40pt
\vsize=54pc
\voffset=6pt

\def\*{{\textstyle *}}
\def\s*{{\scriptstyle *}}
\font\ssanss=plss8 scaled 900
\font\smss=plss8 scaled 800
\def\sss#1{\hbox{\ssanss #1}}
\def\proof{\demo{Proof}}
\def\endproof{\hfill \vrule height4pt width6pt depth2pt \enddemo}
\def\R{{\bb R}}
\def\C{{\bb C}}
\def\r{\ss r}
\def\sT{\ss T}
\def\ssT{\sss T}
\def\smT{\hbox{\smss T}}
\def\ezt{\eu\zt}
\def\ezp{\eu\zp}

\def\xd{\text{\tenrm d}}
\def\sdt{\text{\sevenrm d}_{\smT}}
\def\dt{\xd_{\sss T}}
\def\vt{\text{\tenrm v}_{\sss T}}
\def\iv{\text{\tenrm i}_{\text{\sevenrm v}_{\smT}(\zm)}}
\def\id{\text{\tenrm i}_{\text{\sevenrm d}_{\smT}(\zm)}}

\def\xi{\text{i}}
\def\xit{\xi_{\ssT}}

\def\ll{\text{\it \$}}

\def\blb{\ref \key}
\def\endbb{\endref}
\def\today{\ifcase\month\or
January\or February\or March\or April\or May\or June\or July
\or August\or September\or October\or November\or December\fi
\space\number\day, \number\year}
\newsymbol\ppr 226E
\newsymbol\leqs 1336
\newsymbol\geqs 133E
\topmatter
        \title
        Tangent Lifts of Poisson and Related Structures
        \endtitle

    \author
        J Grabowski$^1$ and
        P Urba\'nski$^2$
    \endauthor
    \affil
     $^1$ Institute of Mathematics, University of Warsaw\\
    Banacha 2, 02-097 Warsaw, Poland\\
     Erwin Schr\"odinger Institute for Mathematical Physics\\
     Pasteurgasse 6/7, 1090 Wien, Austria\\ \\
    $^2$    Division of Mathematical Methods in Physics,
        University of Warsaw \\
        Ho\D{z}a 74, 00-682 Warsaw, Poland
    \endaffil
     \address
     Institute of Mathematics, University of Warsaw, Banacha 2,
     02-097 Warsaw, Poland
     \endaddress
     \email jagrab\@mimuw.edu.pl
     \endemail
     \address
     Division of Mathematical Methods in Physics, University of
     Warsaw, Ho\D{z}a 74, 00-682 Warsaw, Poland
     \endaddress
    \email urbanski\@fuw.edu.pl \endemail

   \thanks
   Classification PACS: 02.40+m, 03.20+i.
   \endthanks
   \thanks
   $^1$ Partially supported by the Isaac Newton Institute for
   Mathematical Sciences, Cambridge, England.
   \endthanks
   \thanks
   $^2$ Supported by KBN, grant No 2 11129 91 01.
   \endthanks
   \abstract
   The derivation $d_T$ on the exterior algebra of forms on a
   manifold $M$ with values in the exterior algebra of forms on the
   tangent bundle $TM$ is extended to multivector fields.
   These tangent lifts are studied with applications to the theory of
   Poisson structures, their symplectic foliations, canonical
   vector fields and Poisson-Lie groups.
   \endabstract
    \endtopmatter
    \leftheadtext{Grabowski and Urba\'nski}
    \rightheadtext{Tangent Lifts}
    \document
    \subhead{0. Introduction}\endsubhead
    A derivation $\dt$ on the exterior algebra of forms on a
manifold $M$ with values in the exterior algebra of forms of the
tangent bundle $\sT M$ plays essential r\^ole in the calculus of
variations (\cite{Tu}) and, in particular, in analytical
mechanics. The derivation $\dt \zw$ of the symplectic 2-form of
a symplectic manifold $(M,\zw)$ provides the tangent bundle $\sT
M$ with a symplectic structure. A vector field $X\colon
M\rightarrow \sT M$ is locally Hamiltonian if its image $X(M)$
is a Lagrangian submanifold of $(\sT M, \dt \zw )$. The concept of a
generalized hamiltonian system can be introduced as a Lagrangian
submanifold of $(\sT M, \dt \zw)$. The infinitesimal  dynamics
of a relativistic particle is an example of such a system. The
derivation $\dt$ has also an aspect of the total   Lie
derivative in the exterior algebra of forms:  $\ll_X\zm =
X^{\textstyle *} \dt\zm $ (Theorem 3.2).

     In analytical mechanics  Poisson structures play the role
as important as symplectic structures. The phase space is
considered as a manifold equipped with a Poisson structure
rather than symplectic one.  On the other hand, in the theory of
systems with symmetries, much attention is paid to the case of
Poisson symmetries, i.~e., the symmetry group is a Poisson-Lie
group.   Poisson-Lie groups are of interest also because of their
relation to quantum groups.

    A Poisson structure is usually given by a bivector field $\zL$
and, in general, not by a two-form and vanishing of the Schouten
bracket corresponds to vanishing of the exterior derivative. This
shows that, in order to generalize the mentioned ideas and
results from the symplectic to the Poisson case, we need to
carry over the discussion from forms to multivector fields.
    The aim of this paper is to extend $\dt$ to the exterior
algebra of multivector fields and to establish relations concerning
Poisson structures which correspond to the mentioned above
relations in symplectic geometry. We would like to emphasize
that our goal is not to extend the general theory of derivations
of forms or vector-valued forms (like in \cite{MCS1-3,By}) to
the case of multivector fields, but to get an analogue of $\dt$
only.

In the first two sections we concentrate  on the definition of $\dt$
on the exterior algebra of forms. The usual definition of $\dt$
as a commutator $[\xit,\xd]$ does not emphasize the r\^ole of the
tangent functor and, since it uses the exterior derivative
$\xd$, cannot be generalized to the case of multivector fields.
An $r$-form $\zm$ on $M$ defines a number of vector bundle morphisms
        $$\widetilde{\zm}^i \colon {\bigwedge}^i \sT M \rightarrow
{\bigwedge}^{r-i} \sT^{\textstyle *}  M.
                                                    $$
    We show that $\widetilde{\dt\zm}^i$ can be obtained from
the tangent morphism $\sT \widetilde{\zm}^i$ by natural
transformation (\cite{KMS}) of functors, which generalize the
well-known natural transformations  $\za_M\colon
\sT\sT^{\textstyle *}M \rightarrow \sT^{\textstyle *} \sT M $
 and $\zk_M\colon \sT\sT M\rightarrow \sT\sT M$.
    With this fact, the generalization of $\dt$ to the case of
multivector fields becomes obvious and it is given in Section~2
(Theorem~2.2). It appears that $\dt$ is a $\vt$-derivation of the
exterior algebra of multivector fields on $M$ with values in the
exterior algebra of multivector fields on $\sT M$, where $\vt$
is the vertical lift. In the case of a vector field $X$ on $M$
the resulting field $\dt X$ is the well-known complete lift (see
e.~g. \cite{MFL}).

 The basic property of $\dt$  is that it commutes with the
Schouten bracket (Theorem~2.5), what corresponds to the fact that $\dt $
commutes with $\xd$ on forms. Thus, if $\zL$ is a bivector field
representing a Poisson structure on $M$, then $\dt\zL$ defines a
Poisson structure on $\sT M$.  We observe that $\dt \zL$ is the
tangent Poisson structure discussed in \cite{SdA,Co1}.
In the presented approach certain functorial properties of $\dt$
become quite obvious.

In Section 3 we show  that, as in
the case of forms, $\dt$ plays the r\^ole of the total Lie
derivative.  As a consequence, we can describe  in Section 6 a
canonical vector field on a Poisson manifold as a Lagrangian
submanifold with respect to the tangent Poisson structure
introduced in Section 5 (compare with \cite{SdA}).
The presentation of the tangent Poisson structure in Section 5 is
close to the one given by Courant in \cite{Co1,Co2,Co3}.

    The derivation $\dt$ helps to identify vector bundle
morphisms $\zn\colon \sT ^{\textstyle *}M \rightarrow \sT M$,
which correspond to  Poisson structures (Theorem 4.4). This
identification is complementary to ones expressed in terms of
the Jacobi identity and of the  Schouten bracket. What is important,
the condition for $\zn$ is expressed in terms of objects and
morphisms and does not require any additional general operations
like exterior derivative and the Schouten bracket. Hence, it is a
subject for functorial treatments.

The remaining part of the paper is devoted to the tangent lift
of Poisson-Lie structures and to the analysis of the symplectic
foliations of tangent Poisson manifolds. We show in Section~7
that the tangent group of a Poisson-Lie group $(G,\zL)$ with the
tangent Poisson structure is again a Poisson-Lie group.  Its
Poisson-Lie algebra is the tangent Poisson-Lie algebra of the
Poisson-Lie algebra of $(G,\zL)$ (Section~8).   In Section~9 we
define the tangent lift of a generalized foliation and in
Section~10 we prove that the symplectic foliation of the tangent
Poisson manifold $(\sT M,\dt\zL)$ is  the tangent lift of the
symplectic foliation of $(M,\zL)$.

    This work is a contribution to a program of geometric
formulations of physical theories conducted jointly with
W.~M.~Tulczyjew.

    The authors wish to thank G.~Marmo for helpful suggestions
concerning the case of Poisson-Lie structures.

\subhead{1. Geometric preliminaries}\endsubhead

    In this section we define morphisms
    $$ \zk^r_M \colon {\bigwedge} ^r \sT\sT M \rightarrow \sT
{\bigwedge} ^r \sT M $$
    and
            $$ \ze^r_M\colon {\bigwedge}^r \sT^{\textstyle *}
\sT M \rightarrow  \sT {\bigwedge}^r \sT^{\textstyle *} M$$
    which generalize the well-known isomorphisms $\zk_M$ and
$\ze_M=\za^{-1}_M$.
Functorial properties of these mappings and their duals are
discussed.

Let $M$  be a smooth manifold. By $\ezt_M \colon \sT M\rightarrow M $ we denote the
tangent fibration and by $\ezp_M\colon \sT^{\textstyle *} M\rightarrow M$ the
cotangent fibration.  For $r=0,1,2,\dots $, we define exterior product bundles
$\bigwedge ^r\sT M$ and $\bigwedge ^r \sT^{\textstyle *}M $ with the canonical
projections $\ezt^r_M \colon \bigwedge ^r \sT M \rightarrow M $ and $\ezp^r_M \colon
\bigwedge ^r \sT^{\textstyle *}  M \rightarrow M $ respectively. For $ r=0$ we have
$\bigwedge ^0 \sT M\simeq \bigwedge ^0 \sT^{\textstyle *} M \simeq M\times \R$.

There is a collection of canonical pairings
        $$  \langle ,\rangle^r_M \colon   {\bigwedge} ^r \sT M \times _{M} {\bigwedge}
^r \sT^{\textstyle *}  M \rightarrow \R .
                                                            $$
By applying the tangent functor to these pairings, we obtain tangent pairings
$$  \langle ,{\rangle'}_{\ssT M}^r \colon   \sT{\bigwedge} ^r \sT M \times _{\ssT M}\sT {\bigwedge}
    ^r \sT^{\textstyle *}  M \rightarrow \sT \R \rightarrow \sT_0\R = \R .
                                                            $$
We used the canonical identification of bundles
        $$ \sT {\bigwedge} ^r \sT M  \times _{\ssT M}\sT{\bigwedge} ^r \sT^{\textstyle *}
M \simeq \sT ({\bigwedge} ^r \sT M  \times _M {\bigwedge} ^r \sT^{\textstyle *}  M ) . $$

Let $(x^i)$ be a local coordinate system in $M$. In bundles $\bigwedge ^r \sT M,\
\bigwedge ^r \sT^{\textstyle *}  M,\ \sT \bigwedge ^r \sT M  $ and $\sT \bigwedge ^r
\sT^{\textstyle *} M$ we have adopted coordinate systems\newline
\centerline{$(x^i, {\dot x}^{j_1\dots j_r})$, $(x^i, { p}_{j_1\dots j_r})$, $(x^i, {\dot
x}^{j_1\dots j_r}, \zd x^k, \zd{\dot x}^{l_1\dots l_r})$ and
  $(x^i, { p}_{j_1\dots j_r}, \dot{ x}^k, {\dot p}_{l_1\dots l_r})$} \newline
respectively, whe\-re $j_1<j_2<\cdots <j_r,$ etc. .
    In these coordinates the introduced pairings read as follows:
        $$ \langle ,\rangle ^r_M \colon ((x^i, {\dot x}^{j_1\dots j_r}),(x^i, { p}_{l_1\dots
l_r}))\mapsto  \sum_{j_1<j_2<\cdots<j_r}{\dot x}^{j_1\dots j_r} p_{j_1\dots j_r}
                                                            $$
and
        $$\multline  \langle ,{\rangle'}_{\ssT M}^r \colon ((x^i, {\dot
x}^{j_1\dots j_r}, \zd x^k, \zd{\dot x}^{l_1\dots l_r}),
  (x^i, { p}_{j_1\dots j_r}, \zd { x}^k, {\dot p}_{l_1\dots l_r})) \mapsto \\
    \sum_{j_1<\cdots <j_r}(\zd\dot x^{j_1\dots j_r}\cdot p_{j_1\dots j_r} +
    \dot x^{j_1\dots j_r}\cdot \dot p_{j_1\dots j_r} ).
                                        \endmultline            $$

    For each manifold $M $ there is a canonical diffeomorphism (cf. \cite{Tu1})
        $$ \zk_M \colon \sT\sT M \rightarrow \sT\sT M
                                                            $$

which is an isomorphism of vector bundles
        $$ \ezt_{\ssT M}\colon \sT\sT M\rightarrow \sT M \ \ \text{and} \ \ \sT\ezt_M
\colon \sT\sT M\rightarrow \sT M  .
                                                            $$
In particular,
        $$\ezt_{\ssT M}\circ \zk_M = \sT\zt_M \ \ \text{and} \
\ \sT\ezt_M\circ \zk_M =\zt_{\ssT M}.
                                                            $$
Regarded as a diffeomorphism of  $\sT\sT M$, \ $\zk_M$ is
involutive: $\zk^2_M = Id_{\ssT\ssT M}$.
By $\za_M$ we denote the isomorphism
        $$ \za_M \colon \sT\sT^{\textstyle *}  M \rightarrow \sT^{\textstyle *} \sT M
                                                                $$
     of vector bundles

        $$  \sT \ezp_{ M}\colon \sT\sT^{\textstyle *}  M\rightarrow \sT M \ \
\text{and} \ \ \ezp_{\ssT M} \colon \sT^{\textstyle *} \sT M\rightarrow \sT M  ,
                                                            $$
 dual to $\zk_M$  with respect to pairings  $\langle ,\rangle '_{\ssT M} = {\langle
,\rangle '}^1_{\ssT M} $  and $\langle ,\rangle _{\ssT M} = \langle ,\rangle
^1_{\ssT M}$:
        $$ \langle v,\za_M (w)\rangle_{\ssT M} = \langle \zk_M(v), w\rangle '_{\ssT M}.
                                                                $$
     In the following, by $\ze_M$ we denote $\za_M^{-1}$.  In the introduced above
local coordinates in $\sT\sT M$, $\sT\sT^{\textstyle *} M$ and the adopted from
$(x^i)$ coordinates $(x^i, \dot x^j, \zp_k, \dot\zp_l)$ in
$\sT^{\textstyle *} \sT M,$ we
have
        $$ (x^i,\dot x^j, \zd x^k, \zd\dot x^l)\circ \zk_M = (x^i, \zd x^j, \dot
x^k, \zd\dot x^l)
                                                                $$
and
        $$(x^i, \dot x^j, \zp_k, \dot \zp_l)\circ \za_M = (x^i, \dot x^j, \dot p_k, p_l).
                                                                $$

Now, we generalize these morphisms to the multilinear case.
Let $\bigwedge ^r_{ M}$ be the wedge product mapping
        $$ {\bigwedge} ^r_{ M} \colon  \times ^r_{\ezt_M} \sT M \rightarrow
{\bigwedge} ^r \sT M .
                                                                $$
We apply the tangent functor to this mapping and we get
        $$ \sT{\bigwedge} ^r_{ M} \colon \sT \times ^r_{\ezt_M} \sT M \rightarrow
\sT{\bigwedge} ^r \sT M .
                                                                $$
Since $\zk_M$ extends to an isomorphism of vector bundles
        $$ \times ^r\zk_M \colon  \times ^r_{\ezt_{\smT M}} \sT\sT M \rightarrow
\times ^r_{\ssT \ezt_M}\sT\sT M \simeq \sT \times ^r_{\ezt_M} \sT M,
                                                                $$
we get also
        $$ \sT{\bigwedge} ^r_{M} \circ (\times ^r{\zk_M})\colon \times
^r_{\ezt_{\smT M}}\sT\sT M\rightarrow  \sT {\bigwedge} ^r \sT M .
                                                                $$

It is easy to verify that this mapping is multilinear and skew-symmetric and,
consequently, defines a  morphism  $\zk^r_M$ of vector bundles over $\sT M$:
        $$ \zk^r_M \colon {\bigwedge} ^r \sT\sT M \rightarrow \sT {\bigwedge} ^r \sT M.
                                                                \tag 1.1$$
In other words,
        $$ \zk^r_M \circ {\bigwedge} ^r_{\ssT M} = \sT{\bigwedge}^r_M \circ
\times^r\zk_M ,
                                                                $$
i. e., the following diagram is commutative

        $$ \CD
   \times ^r_{\ezt_{\smT M}}\sT \sT M @> \dsize \wedge^r_{\ssT M} >> \bigwedge^r \sT\sT M \\
   @V\dsize \times ^r\zk_M VV                 @VV\dsize \zk^r_M V \\
   \sT\times ^r_{\ezt_M}\sT M   @>\dsize \sT\wedge ^r_M >>  \sT\bigwedge^r \sT\sT M
                                                \endCD \quad .
                                                        \tag 1.2$$

Of course, $\zk^1_M = \zk_M$ and
for a simple $r$-vector $v_1\wedge \dots \wedge v_r$ on $\sT M$ we have
        $$ \zk^r_M (v_1\wedge \dots \wedge v_r) = \sT{\bigwedge} ^r_{M}
(\zk_M(v_1),\dots ,\zk_M(v^r)).
                                                                $$

In local coordinates $\zk^r_M $ reads as follows:
        $$ (x^i, \dot x^{j_1\dots j_r}, \zd x^k, \zd\dot x^{l_1\dots l_r})\circ
\zk^r_M = (x_i, \zd x^{j_1\dots j_r}, \dot x^k, \sum_m \zd x^{l_1\dots l_{m-1}}\wedge
\zd\dot x^{l_m} \wedge \zd x^{l_{m+1}\dots l_r }),
                                                                $$
where $(x^i,\dot x^j, \zd x^{k_1\dots k_n}\wedge \zd\dot x^{l_{n+1}\dots l_r})$ are
adopted coordinates in $\wedge ^r \sT\sT M$, with the obvious identification
        $$\zd x^{l_1\dots l_{m-1}}\wedge \zd\dot x^{l_m} \wedge \zd x^{ l_{m+1}
\dots l_r } = (-1)^{r-m}\zd x^{l_1\dots l_{m-1}l_{m+1}\dots l_r}\wedge \zd\dot x^{l_m}.
                                                                $$

Let $\phantom{\bigwedge}^i{\bigwedge}^{r-i}_M $ be the wedge product
        $$ \phantom{\bigwedge}^i{\bigwedge}^{r-i}_M \colon {\bigwedge}^i\sT
M\times _M  {\bigwedge}^{r-i}\sT M \rightarrow {\bigwedge}^r\sT M .$$
 From the diagram (1.2) we easily get that the diagram

        $$ \CD {\bigwedge}^i \sT\sT M \times _{\ssT M} {\bigwedge}^{r-i} \sT\sT M
@>\dsize \zk^i_M \times \zk^{r-i}_M >> \sT{\bigwedge}^i \sT M \times _{\ssT M}\sT
{\bigwedge}^{r-i} \sT M \\ @V\dsize \vphantom{\wedge}^i{\wedge}^{r-i}_{\ssT M} VV
@V\dsize \sT \vphantom{\wedge}^i{\wedge}^{r-i}_M VV \\ {\bigwedge}^r \sT\sT M
@>\dsize \hphantom{aaaa} \zk^r_M \hphantom{aaaa} >> \sT{\bigwedge}^r\sT M
                                                \endCD  \tag 1.3    $$

is commutative.

In order to define $\ze^r_M$, a counterpart to $\zk^r_M$, we consider first
the mapping
        $$\sT {\bigwedge}^r_{\ssT ^\s* M} \circ (\times^r\ze_M) \colon \times
^r_{\ezp_{\ssT M}} \sT ^{\textstyle *} \sT M \rightarrow \sT {\bigwedge}^r
\sT^{\textstyle *} \sT M.
                                                                $$
Since it is multilinear and skew-symmetric, it defines a mapping
        $$ \ze^r_M\colon {\bigwedge}^r \sT^{\textstyle *} \sT M \rightarrow  \sT
{\bigwedge}^r \sT^{\textstyle *} M
                                                                $$
and, as in the case of $\zk^r_M$, we have
        $$\ze^r_M \circ {\bigwedge}^r_{\ezp_{\ssT M}} = \sT {\bigwedge}^r_{\ezp_M}
\circ  \times ^r \ze_M,
                                                                $$
     i. e., the diagram
        $$ \CD   \times ^r_{\ezp_{\smT M}}\sT^{\textstyle *}  \sT M @>\dsize
\wedge^r_{\ssT M} >> \bigwedge^r \sT ^{\textstyle *} \sT M \\
   @V\dsize \times ^r\ze_M VV                 @VV \dsize\ze^r_M V \\
   \sT\times ^r_{\ezt_M}\sT^{\textstyle *}  M   @>\dsize \sT\wedge ^r_M >>
\sT\bigwedge^r \sT ^{\textstyle *}   M
                                            \endCD          \tag 1.4    $$
is commutative.

Let $(x^i, \dot x^j, \zp^{k_1\dots k_n}\wedge \dot\zp ^{l_{n+1}\dots l_r} )$ be the
adopted coordinate system in ${\bigwedge}^r \sT^{\textstyle *} \sT M$. We have
        $$(x^i, p_{j_1\dots j_r}, \dot x^k, \dot p_{l_1\dots l_r})\circ \ze^r_M
 = (x^i, \dot\zp_{j_1\dots j_r}, \dot x^k, \sum_m \dot\zp_{l_1\dots l_{m-1}}\wedge
\zp_{l_m}\wedge \dot\zp_{l_{m+1}\dots l_r} ).
                                                                $$
Here we identified $\dot\zp_{l_1\dots \l_{m-1}}\wedge \zp_{l_m}\wedge
\dot\zp_{l_{m+1}\dots l_r}$ and $(-1)^{m-1}\zp_{l_m} \wedge \dot\zp_{l_1\dots
l_{m-1} l_{m+1} \dots l_r} $.

We  have also a commutative diagram
        $$ \CD {\bigwedge}^i \sT^{\textstyle *} \sT M \times _{\ssT M}
{\bigwedge}^{r-i} \sT^{\textstyle *} \sT M  @>\dsize \ze^i_M \times \ze^{r-i}_M >>
\sT{\bigwedge}^i \sT^{\textstyle *}  M \times _{\ssT M}\sT {\bigwedge}^{r-i}
\sT^{\textstyle *}  M \\ @V\dsize \phantom{\wedge}^i{\wedge}^{r-i}_{\ssT M} VV  @A
{\dsize \sT \vphantom{\wedge}^i{\wedge}^{r-i}_M} AA \\ {\bigwedge}^r \sT^{\textstyle
*} \sT M @>\dsize \hphantom{aaaa} \ze^r_M \hphantom{aaaa} >>
\sT{\bigwedge}^r\sT^{\textstyle *}  M
                                            \endCD \quad .          \tag 1.5    $$

By ${\zk^r_M}'$ and ${\ze^r_M}'$ we denote  vector bundle morphism
        $$ {\zk^r_M}' \colon  \sT {\bigwedge}^r \sT^{\textstyle *}   M
\rightarrow  {\bigwedge}^r \sT^{\textstyle *}  \sT M,
                                                                $$
and
        $$ {\ze^r_M}' \colon  \sT {\bigwedge}^r \sT   M
\rightarrow  {\bigwedge}^r \sT \sT M,
                                                                $$
dual  to $\zk^r_M$  and $\ze^r_M$ with respect to pairings $\langle ,{\rangle'_{\ssT M}}^r$ and
$\langle ,\rangle_{\ssT M}^r$.
 We have, in particular, ${\zk^1_M}' = \zk_M' = \ze_M^{-1} = \za_M$ and ${\ze^1_M}'
= \zk_M$.

\specialhead Functorial properties \endspecialhead
It is known that $\zk_M$ and $\za_M$ are natural transformations of iterated functors
$\sT\sT,\ \sT\sT^{\textstyle *}$ and $\sT^{\textstyle *} \sT $, i.~e., that for every
morphism $ \zf\colon M\rightarrow N$ we have
        $$ \split \zk_N\circ \sT\sT\zf &= \sT\sT\zf \circ \zk_M, \\
    \za_M\circ \sT\sT^{\textstyle *} \zf &= \sT^{\textstyle *}
     \sT\zf \circ \za_N, \\
    \sT\sT^{\textstyle *}  \zf\circ \ze_N &=
     \ze_M\circ \sT^{\textstyle *} \sT\zf .
                                            \endsplit\tag 1.6       $$
Note that    $\sT^{\textstyle *} \zf$ is, in general, not a mapping but a relation only,
with the domain $\sT^{\textstyle *} _{\zf(M)}N$ and codomain $(\ker \sT\zf)^o$.
    Morphisms
        $$ {\bigwedge}^r\sT \zf \colon {\bigwedge}^r\sT M\rightarrow {\bigwedge}^r
\sT N
                                                                $$
and
        $$ {\bigwedge}^r\sT^{\textstyle *}  \zf \colon {\bigwedge}^r
\sT^{\textstyle *}  N\rightarrow {\bigwedge}^r \sT^{\textstyle *}  M
                                                                $$
are defined by relations
        $$ {\bigwedge}^r\sT\zf \circ {\bigwedge}^r_M = {\bigwedge}^r_N \circ
\times ^r \sT \zf
                                                                $$
and
        $$ {\bigwedge}^r\sT ^{\textstyle *} \zf \circ {\bigwedge}^r_N =
{\bigwedge}^r_M \circ \times ^r \sT ^{\textstyle *}  \zf .
                                                                $$

\proclaim{Theorem 1.1}{}
For $\zf \colon  M\rightarrow N$ we have
        $$ \zk^r_N \circ {\bigwedge}^r \sT\sT \zf =  \sT {\bigwedge}^r\sT\zf \circ
\zk^r_M
                                                                $$
and
        $$\sT {\bigwedge}^r \sT^{\textstyle *} \zf\circ \ze^r_N = \ze^r_M \circ
{\bigwedge}^r \sT^{\textstyle *} \sT\zf .
                                                                $$
\endproclaim
\proof
 From the definition of $\zk^r_N$ it follows that
        $$\multline \zk^r_N \circ {\bigwedge}^r \sT \sT\zf  \circ
{\bigwedge}^r_{\ssT M} = \zk^r_N\circ  {\bigwedge}^r_{\ssT N}\circ  \times ^r\sT \sT
\zf = \sT {\bigwedge}^r_{ N}\circ \times ^r\zk_N \circ \times ^r \sT\sT\zf \\
    =\sT {\bigwedge}^r_{ N}\circ \times ^r(\zk_N \circ \sT\sT \zf) =
\sT{\bigwedge}^r_{ N} \circ \times ^r (\sT\sT  \zf \circ \zk_M) =
\sT{\bigwedge}^r_{N}\circ \times ^r \sT\sT \zf \circ \times ^r\zk_M .
                                    \endmultline                    $$

Since
        $$ {\bigwedge}^r \sT \zf \circ {\bigwedge}^r_M = {\bigwedge}^r_N\circ
\times ^r \sT \zf,
                                                                $$
we get
        $$\sT{\bigwedge}^r \sT \zf \circ \sT {\bigwedge}^r_{ M} = \sT
{\bigwedge}^r_{ N} \circ \times ^r \sT\sT\zf
                                                                $$
and
        $$ \sT{\bigwedge}^r_{N} \circ \times ^r \sT\sT\zf \circ \times ^r\zk_M =
\sT {\bigwedge}^r \sT \zf
\circ \sT{\bigwedge}^r_{M} \circ \times ^r\zk_M = \sT {\bigwedge}^r \sT \zf \circ
\zk^r_M\circ  {\bigwedge}^r_{\ssT M}.
                                                                $$
This completes the proof of the first identity. The proof of the second one is analogous.
\endproof

We have also the dual identities:
\proclaim{Theorem 1.2}{}
For $\zf \colon  M\rightarrow N$ we have
        $$ {\zk^r_M}' \circ\sT {\bigwedge}^r \sT^{\textstyle *}  \zf =
{\bigwedge}^r \sT^{\textstyle *} \sT\zf \circ  {\zk^r_N}'
                                                                $$
and
        $$ {\bigwedge}^r \sT \sT \zf\circ {\ze^r_M}' = {\ze^r_N}' \circ
\sT {\bigwedge}^r \sT\zf .
                                                                $$
\endproclaim

\subhead 2. Derivation $\dt$ of differential forms and multivector fields \endsubhead

In this section we refer to the theory of derivations of
differential forms as presented in~\cite{PiTu}. We define the
derivation  $\dt$ on forms in a way which differs from the standard one,
but which shows its obvious extension to the case of multivector
fields. It appears that defined operation  $\dt$ on multivector
fields is a derivation of degree zero with respect to the
vertical lift of multivector fields. The most important property
of $\dt$ is that it commutes with the Schouten bracket.
    \proclaim{Definition}{} Let $\zF=\bigoplus^{\infty}_{q=0}\zF^q$ and
$\Psi=\bigoplus^{\infty}_{q=0}\Psi^q $ be commutative graded algebras and let
$\zr\colon \zF \rightarrow \Psi$ be a graded algebra
homomorphism. A linear mapping $a\colon \zF\rightarrow \Psi$ is
called a $\zr$-{\rm derivation} of degree $r$ if:
        $$a(\zF^q)\subset \Psi^{q+r} $$
and
        $$a(\zm \wedge \zn) =a(\zm)\wedge\zr( \zn) +
(-1)^{qr}\zr(\zm)\wedge a(\zn), $$
where $q=$ degree $\zm$.
    \endproclaim
Let $M$ be a manifold and
let $\zt \colon E \rightarrow M$ be a vector fibration. By
$\zF(\zt)$ we denote the graded  exterior algebra generated by
sections of $\zt$. For $\zt = \ezp_M$ we get the graded algebra
of forms on the manifold $M$ and for $\zt = \ezt_M$- the graded
algebra of multivector fields on $M$.

Let be $\zm \in \zF^r(\ezp_M)$, i.~e., $\zm$ is an $r$-form on $M$. The {\it vertical
lift} of $\zm$ is an $r$-form  $\vt (\zm)\in \zF^r(\ezp_{\ssT M})$, $\vt(\zm) =
\ezt^{\textstyle *} _M \zm$, i.~e., $\vt(\zm) $ is the pull-back of $\zm$ with respect to
the projection $\ezt_M$. Since the pull-back commutes with the wedge product, the mapping
        $$ \vt\colon \zF(\ezp_M) \rightarrow \zF(\ezp_{\ssT M})
                                                                $$
      is a homomorphism of graded commutative algebras.

A {\it second order vector field } $\zG$ on $M$ is a vector field on $\sT M$  such
that
        $$ \ezt_{\ssT M}\circ \zG = \text{\tenrm id}_{\ssT M}= (\sT\zt_M)\circ \zG
                                                                $$
or, equivalently, $\zk_M \circ \zG = \zG$.
In adopted local coordinates,
        $$ \zG(x,\dot x) = \dot x^k \frac{\partial}{\partial x^k} + f^k(x,\dot x)
\frac{\partial}{\partial \dot x^k}.
                                                                $$
The contraction of the vertical lift of  a form $\zm \in \zF^r(\ezp_M)$, $r>0$, with $\zG$
does not depend on the choice of the second order field $\zG$. We define
        $$\xit {\zm} =\cases \xi_\zG (\vt\zm) \in \zF^{r-1}(\ezp_{\ssT M}),
&\text{for $r>0$}\\
                        0 &\text{for $r=0$}. \endcases
                                                                $$

The {\it tangent lift } $\dt \zm$ of $\zm \in \zF(\ezp_M)$ is defined by
        $$ \dt \zm = \xd \xit \zm + \xit \xd \zm = \ll_\zG\vt(\zm).
                                                                $$

Since $\ll_\zG$ is a derivation in $\zF(\ezp_{\ssT M})$ and $ \vt$ is a homomorphism
of graded algebras, it follows that $\dt$ is a $\vt$-derivation of degree $0$.

If, in local coordinates, $\zm = \zm_{i_1\dots i_r}\xd x^{i_1}\wedge
\cdots\wedge \xd x^{i_r},$ then
        $$ \dt\zm (x,\dot x) = \frac{\partial \zm_{i_1\dots i_r}}{\partial x^k}(x)
  \dot x^k \xd x^{i_1}\wedge \dots \wedge \xd x^{i_r} + \sum_m
\zm_{i_1\dots i_r}(x) \xd x^{i_1}\wedge \dots\wedge \xd \dot x^{i_m}\wedge \dots \xd x^{i_r}
                                                        \tag 2.1    $$
for $r>0$ and
        $$ \dt\zm(x, \dot x) = \frac{\partial \zm}{\partial x^i}(x)\dot x^i
                                                                $$
for $r=0$.

The operation $\xit$, which is, in fact, a $\vt$-derivation of degree $-1$, can
be defined in a more intrinsic way. An $r$-form $\zm$, $r>0$, defines a vector bundle morphism
        $$ \widetilde{\zm}^1\colon \sT M \rightarrow {\bigwedge}^{r-1}
\sT^{\textstyle *}  M \colon v\mapsto \xi_{\textstyle v}\zm
                                                                $$
     and the following formula holds:
        $$ \xit \zm = (\widetilde{\zm}^1)^{\textstyle *} \zvy^{r-1}_M,
                                                        \tag 2.2    $$
where $\zvy^{r-1}_M $ is the canonical (Liouville) $(r-1)$-form on ${\bigwedge}^{r-1}
\sT^{\textstyle *}  M$. The Liouville form is defined by
        $$ \zvy_M^{r-1}(a)(v_1,\dots, v_{r-1}) = a(\sT\ezp_M(v_1),\dots
\sT\ezp_M(v_{r-1})) .
                                                                $$

Let us notice also that for $r=0$ we have $\dt\zm(v)= \langle v,\xd \zm\rangle $
($\zm$ is a function).

The tangent lift $\dt \zm$ can be defined more directly by means of the tangent
functor. Let us fix $0 \leqs i \leqs r$. An $r$-form $\zm$ on $M$ defines, in an
obvious way, a vector bundle morphism
        $$\widetilde{\zm}^i \colon {\bigwedge}^i \sT M \rightarrow
{\bigwedge}^{r-i} \sT^{\textstyle *}  M.
                                                                $$

Now, we define $\zk^r_M$ and $ \ze^r_M$ for $r=0$. We have ${\bigwedge}^0\sT M =
{\bigwedge}^0 \sT^{\textstyle *}  M = M\times \R$.
We define
        $$ \zk^0_M = \ze^0_M \colon {\bigwedge}^0 \sT\sT M =  \sT M\times \R
\rightarrow  \sT {\bigwedge}^0 \sT M = \sT(M\times \R) = \sT M \times \R\times \sT_0\R
                                                                $$
by
        $$ (x^i,\dot x^j, t,\dot t)\circ \zk^0_M = (x^i, \dot x^j, 0,t).
                                                                $$
The dual mapping ${\zk^0_M}' \colon \sT{\bigwedge}^0 \sT M \rightarrow
{\bigwedge}^0 \sT\sT M $ is given by
        $$ (x^i,\dot x^j, t)\circ {\zk_M^0}' = (x^i,\dot x^j, \dot t).
                                                                $$

    \proclaim{Theorem 2.1}{}
    Let be $\zm \in \zF^r(\ezp_M)$ and  $0 \leqs i\leqs r $. The following diagram is
commutative
        $$ \CD  \sT{\bigwedge}^i \sT M @>\dsize\sT\widetilde{\zm}^i >>
\sT{\bigwedge}^{r-i} \sT^{\textstyle *}  M\\
    @A\dsize\zk^i_M AA     @V\dsize (\zk_M^{r-i})'VV \\
    {\bigwedge}^i \sT \sT M @>\dsize\widetilde{\dt\zm}^i >> {\bigwedge}^{r-i}
\sT^{\textstyle *} \sT M
                                            \endCD \ \ .    \tag 2.3$$

    \endproclaim
\proof
We show first that there exists $D(\zm)\in \zF^r(\ezp_{\ssT M})$ such that
        $$\widetilde{D({\zm})}^i = (\zk^{r-i}_M)'\circ \sT \widetilde{\zm}^i \circ \zk^i_M.
                                                                $$
     In order to do this, we apply the tangent functor to the commutative diagram

        $$ \CD {\bigwedge}^r \sT M @>\dsize \hphantom{aaaa} \widetilde{\zm}^r
\hphantom{aaaa} >> \R \\
    @A\dsize \phantom{\wedge}^i\wedge^{r-i}_M AA @A\dsize \langle ,\rangle ^{r-i}_M AA \\
    {\bigwedge}^i \sT M \times _M {\bigwedge}^{r-i} \sT M  @> \dsize
\widetilde{\zm}^i\times \text{\tenrm id} >>  {\bigwedge}^{r-i} \sT^{\textstyle *}  M
\times _M {\bigwedge}^{r-i} \sT M
                                                \endCD  \tag 2.4$$

     and, since the diagram (1.5) is commutative, we get the following commutative diagram
        $$\centerline{$ \CD  {\bigwedge}^r \sT \sT M @>\dsize \zk^r_M >> \sT {\bigwedge}^r \sT M @>
\dsize\dt\widetilde{\zm}^r >> \R \\
    @A\dsize \phantom{\wedge}^i\wedge^{r-i}_{\ssT M} AA @A\dsize \sT
\vphantom{\wedge}^i\wedge^{r-i}_M AA @A \dsize\langle ,{\rangle'}_{\ssT M}^{r-i}
A\dsize \qquad (2.5).A\\
{\bigwedge}^i\sT \sT M \times _{{}_{\smT M}}{\bigwedge}^{r-i}\sT \sT M @>\dsize \zk^i_M\times
\zk^{r-i}_M >> \sT{\bigwedge}^i \sT M \times _{{}_{\smT M}} \sT {\bigwedge}^{r-i} \sT M @>
\dsize\sT\widetilde{\zm}^i \times \text{\tenrm id}>> \sT{\bigwedge}^{r-i} \sT^{\textstyle
*}  M \times _{{}_{\smT M}} \sT {\bigwedge}^{r-i} \sT M
                                                \endCD $}
                                                        $$

Here we regard 0-forms as functions rather than sections of
${\bigwedge}^0\sT$--bundles.
It shows that for $u\in {\bigwedge}^i \sT\sT M$ and $v\in {\bigwedge}^{r-i} \sT \sT
M$ we have
        $$ \langle \sT\widetilde{\zm}^i \circ \zk^i_M(u) , \zk^{r-i}_M(v)
{\rangle'}_{\ssT M}^{r-i} = \dt \widetilde{\zm}^r\circ \zk^r_M(u\wedge v)
                                                                $$
    and, consequently,
        $$ \dt \widetilde{\zm}^r\circ \zk^r_M(u\wedge v) = \langle
(\zk_{M}^{r-i})'\circ \sT \widetilde{\zm}^i \circ \zk^i_M(u), v \rangle_{\ssT
M}^{r-i}.
                                                                $$
It follows that $D(\zm) $ exists and $D(\zm)(\cdot ) = \dt \widetilde{\zm}^r\circ
\zk^r_M$.
    Since in local coordinates
        $$ \widetilde{\zm}^r(x,\dot x) = \zm_{i_1\dots i_r}\dot x^{i_1\dots i_r},
                                                                $$
we get
        $$ \dt\widetilde{\zm}^r(x, \dot x, \zd x, \zd\dot x) = \frac{\partial
\zm_{i_1\dots i_r}}{\partial x^k}(x)\zd x^k \dot x^{i_1\dots i_r}  +
\zm_{i_1\dots i_r}(x) \zd \dot x^{i_1\dots i_{r}}
                                                                $$
and
        $$ \dt \widetilde{\zm}^r\circ \zk^r_M =  \frac{\partial
\zm_{i_1\dots i_r}}{\partial x^k}(x)\dot x^k \zd x^{i_1\dots i_r}  + \sum_m
\zm_{i_1\dots i_r}(x) \zd x^{i_1\dots i_{m-1}}\wedge \zd\dot
x^{i_m}\wedge \zd x^{i_{m+1}\dots i_r}.
                                                                $$
The right hand side in this formula corresponds precisely the right hand side in (2.1).
It follows that $\dt \widetilde{\zm}^r\circ \zk^r_M = \widetilde{\dt \zm}^r$ and this
completes the proof.
\endproof

\specialhead Derivations of multivector fields and the Schouten bracket. \endspecialhead
A similar construction can be done in the case of multivector fields. Let be $X\in
\zF^r(\ezt_M)$. As in the case of forms, we have a family of contraction mappings
        $$\widetilde{X}^i \colon {\bigwedge}^i \sT^{\textstyle *}  M \rightarrow
{\bigwedge}^{r-i} \sT M.
                                                                $$

\proclaim{Theorem 2.2}{}
There is a uniquely defined multivector field $\dt X \in \zF^r(\ezt_{\ssT M})$ such
that
        $$ \widetilde{\dt X}^i = (\ze^{r-i}_M)'\circ \sT\widetilde{X}^i \circ
\ze^i_M
                                                                $$
for $i=0,1,\dots, r$, i.~e., the following diagram is commutative

        $$ \CD  \sT{\bigwedge}^i \sT^{\textstyle *}  M @>\dsize \sT\widetilde{X}^i >>
\sT{\bigwedge}^{r-i} \sT  M\\
    @A\dsize \ze^i_M AA     @V\dsize (\ze_M^{r-i})'VV \\
    {\bigwedge}^i \sT^{\textstyle *}  \sT M @>\dsize\widetilde{\dt X}^i >>
{\bigwedge}^{r-i} \sT \sT M
                                                \endCD \quad .      \tag 2.6$$
\endproclaim
The proof goes on like the proof of Theorem 2.1. In particular, we get
        $$\widetilde{\dt X}^r = \dt\widetilde{X}^r \circ \ze^r_M,
                                                                $$
i.~e.,
        $$\dt X(a_1,\dots, a_r) = \dt\widetilde{X}^r( \ze_Ma_1\wedge \dots \wedge
\ze_M a_r).
                                                                $$
Now, writing  in local coordinates
        $$X = X^{i_1\dots i_r}(x)\frac{\partial }{\partial x_{i_1}}\wedge \cdots
\wedge  \frac{\partial }{\partial x_{i_r}}
                                                                $$
and $p = p_{i_1\dots i_r}\xd x^{i_1}\wedge \cdots \wedge \xd x^{i_r}$,we have
        $$\widetilde{X}^r(p) = X^{i_1\dots i_r}p_{i_1\dots i_r}
                                                                $$
and
        $$ \multline \dt\widetilde{X}^r(x^i, \dot x^k, \zp_{i_1\dots i_m}\wedge \dot\zp
_{i_{m+1}\dots i_r}) = \frac{\partial X^{i_1\dots i_r} }{\partial x^{i_k}} (x)\dot
x^k \dot\zp_{i_1\dots i_r} +\\ + \sum_m X^{i_1\dots i_r}(x)\dot\zp_{i_1\dots
i_{m-1}}\wedge \zp_{i_m}\wedge \dot\zp_{i_{m+1}\dots i_r},\endmultline
                                                                $$
i.~e.,
        $$\dt X = \frac{\partial X^{i_1\dots i_r} }{\partial x^{i_k}} (x)\dot
x^k \frac{\partial }{\partial \dot x^{i_1}}\wedge \cdots \wedge \frac{\partial
}{\partial \dot x^{i_r}} + \sum_m X^{i_1\dots i_r}(x) \frac{\partial }{\partial\dot
x^{i_1}} \wedge \cdots \wedge \frac{\partial }{\partial x^{i_m}}\wedge \cdots \wedge
\frac{\partial }{\partial\dot x^{i_r}}.
                                                                $$

Now, let  $X $ be a simple $r$-vector field, i.~e., $X = X_1\wedge \cdots \wedge
X_r$ for some vector fields $X_i\in \zF(\ezt_M)$.
We can consider $\widetilde{X}^r$ as a multilinear, skewsymmetric function on
$\times ^r_M \sT^{\textstyle *} M$
        $$ \widetilde{X}^r(a_1,\dots,a_r) = \sum_{\zs\in S_n} (-1)^\zs
\widetilde{X_{\zs(1)}}^1(a_{1})\cdots \widetilde{X_{\zs(r)}}^1(a_{r}).
                                                                $$
Since  $\dt$ is  a $\vt = \ezt^{\textstyle *} _{\times ^r\ssT ^\s* M}$-derivation on
forms, we have
        $$
\dt\widetilde{X}^r = \sum_{\zs\in S_n} (-1)^\zs \sum_m
\vt(\widetilde{X_{\zs(1)}}^1)\cdots  \dt(\widetilde{X_{\zs(m)}}^1) \cdots \vt(
\widetilde{X_{\zs(r)}}^1)
                                                                $$
      and
        $$ \widetilde{\dt X}^r= \sum_{\zs\in S_n} (-1)^\zs \sum_m
(\vt(\widetilde{X_{\zs(1)}}^1))\circ \ze_M \cdots (\dt( \widetilde{X_{\zs(m)}}^1))
\circ \ze_M \cdots (\vt( \widetilde{X_{\zs(r)}}^1)) \circ \ze_M .
                                                                $$
    Let $Y\in\zF^1(\ezt_M)$ be a vector field on $M$.  Functions
$\vt(\widetilde{Y}^1)\circ \ze_M $ and $\dt(\widetilde{Y}^1)\circ \ze_M$ are  linear
functions on $\sT^{\textstyle *} \sT M$ and define  vector fields on $\sT M$.
In local coordinates we have for $Y = Y^i\dfrac{\partial }{\partial x_{i}}$:
        $$\widetilde{Y}^1(x,p) = Y^i(x)p_i  \ \ \text{and}\ \
\vt(\widetilde{Y}^1)(x,p,\dot x,\dot p) = Y^i(x)p_i.
                                                                $$
Hence, from the definition of $\ze_M$,
        $$ \vt(\widetilde{Y}^1)\circ \ze_M (x,\dot x, \zp, \dot\zp)
          = Y^i(x)\dot\zp_i.
                                                                $$
It follows that $\vt(\widetilde{Y}^1)\circ \ze_M = \widetilde{Y^v}^1$, where $Y^v =
Y^i(x)\dfrac{\partial }{\partial \dot x^{i}}$ is the {\it vertical lift } of $Y$.
The vertical lift $\zF^1(\ezt_M)\ni Y\mapsto Y^v \in \zF^1(\ezt_{\ssT M})$ is linear
and can be extended in a unique way to a homomorphism $\vt$ of graded algebras
        $$\vt \colon \zF(\ezt_M)\rightarrow \zF(\ezt_{\ssT M}).
                                                                $$
The vertical lift $Y^v$ is the generator of a one-parameter group $(\zc^t)$ of
diffeomorphisms (a flow) of $\sT M$ defined by
        $$ \zc^t(v) = v +tY(\ezt_M(v)).
                                                                $$

    In a similar way we get
        $$\dt \widetilde{Y}^1(x,p,\dot x, \dot p) = Y^i(x)\dot p_i +
\frac{\partial Y^i }{\partial x_{k}}(x)p_i\dot x^k
                                                                $$
and
        $$ \dt\widetilde{Y}^1\circ \ze_M (x,\dot x, \zp, \dot\zp) = Y^i(x)\zp_i +
\frac{\partial Y^i}{\partial x^{k}}(x)\dot x^k \dot \zp_i.
                                                                $$
It follows that  $\dt\widetilde{Y}^1\circ \ze_M = \widetilde{Y^c}^1$, where
        $$ Y^c(x,\dot x) = Y^i(x)\frac{\partial }{\partial x^{i}} + \frac{\partial
Y^i}{\partial x^{k}}(x)\dot x^k \frac{\partial }{\partial \dot x^{i}}
                                                                $$
is the {\it complete lift} of $Y$.
The vector field $Y^c$ is the generator of the (local) one-parameter
group of diffeomorphisms $\sT\zf^t\colon \sT M\rightarrow \sT M$, where
$(\zf^t)$ is the flow generated by the vector
field $Y$.

Thus we have proved the following theorem.
    \proclaim{Theorem 2.3}{} The mapping
        $$\dt \colon \zF(\ezt_M)\rightarrow \zF(\ezt_{\ssT M})
                                                                $$
is a $\vt$-derivation of degree $0$ with $\dt (Y)=Y^c$ for
$Y\in\zF ^1(\ezt_M).$

    \endproclaim
The following is  a well-known theorem (e.~g. \cite{MFL,Co1,Co2}, ).
    \proclaim{Theorem 2.4}{}
    The Lie bracket of vertical and complete lifts satisfy the following
commutation relations
        $$ \split [X^v,Y^v] &= 0 ,\\
                [X^c,Y^c] &= [X,Y]^c  ,\\
                [X^v,Y^c] &= [X,Y]^v .\endsplit
                                                        \tag 2.7    $$
    \endproclaim

The Lie bracket of vector fields can be extended to a graded Lie bracket on the
graded space $\zF(\ezt_M)$ of multivector fields -- the Schouten bracket $[\cdot
,\cdot ]$. In this graded Lie algebra the space $\zF^r(\ezt_M)$ is of degree $(r-1).$
Let be $X\in \zF^r(\ezt_M)$. Then
        $$ [X,Y\wedge Z] = [X,Y]\wedge Z + (-1)^{s(r-1)}Y\wedge [X,Z]
                                                                $$
for $Y\in \zF^s(\ezt_M)$, i.~e., $\text{ad}_X $ is a graded derivation of degree
$(r-1)$ of the graded commutative algebra $\zF(\ezt_M)$. The mapping
        $$ \text{ad}\colon (\zF(\ezt_M), [\cdot ,\cdot ]) \rightarrow
\text{Der}(\zF(\ezt_M),\wedge )
                                                                $$
     is a homomorphism of graded algebras.
For simple multivectors we have the following formula (cf. \cite{Mi}):
        $$ [X_1\wedge \cdots\wedge X_p,\, Y_1\wedge  \cdots \wedge Y_q ] = \sum_{i,j}
(-1)^{i+j} [X_i,Y_j]\wedge X_1\wedge \cdots X_{i-1}\wedge X_{i+1}\wedge \cdots
\wedge Y_{j-1}\wedge Y_{j+1}\wedge \cdots Y_q .
                                                        \tag 2.8   $$

    \proclaim{Theorem 2.5}{}
    The derivation $\dt$ commutes with the Schouten bracket, i.~e.,
        $$[\dt X , \dt Y] = \dt [X,Y].
                                                                $$
    \endproclaim
    \proof
Let be $X = X_1\wedge \cdots \wedge X_p$. By $(X_1,\dots,X_p)_n$ we denote a
$p$-vector field
$X_1^v\wedge \cdots \wedge X_n^c\wedge \cdots \wedge X^v_p$ and by
$(X_1,\dots,X_p)_n^i$ a $(p-1)$-vector field as above with the $i$-th factor omitted.
We have then $\dt X = \sum_{i=1}^p (X_1,\dots,X_p)_i$. From the formula (2.8) and
Theorem 2.3 we get, for $Y =Y_1\wedge \dots\wedge Y_q$,
        $$ \multline [\dt X, \dt Y] = \left[ \sum_{n=1}^p (X_1,\dots,X_p)_n ,
\sum_{m=1}^q (Y_1,\dots,Y_p)_m\right] \\= \sum_{n=1}^p\sum_{m=1}^q[ (X_1,\dots,X_p)_n,
(Y_1,\dots, Y_p)_m] = \\
    =\sum_{n=1}^p\sum_{m=1}^q \left( \sum_{j\neq m}(-1)^{n+j} [ X^c_n,Y^v_j]\wedge
(X_1,\dots,X_p)_n^n \wedge (Y_1,\dots,Y_q)^j_m + \right. \\ +
(-1)^{n+m}[X^c_n,Y^c_m] \wedge (X_1,\dots,X_p)_n^n \wedge (Y_1,\dots,Y_p)_m^m +\\
\left. + \sum_{i\neq m}(-1)^{m+i} [ X^v_i,Y^c_m]\wedge (X_1,\dots,X_p)_n^i \wedge
(Y_1,\dots,Y_q)^m_m \right) = \\
    = \sum_{n,j}(-1)^{n+j} [X_n,Y_j]^v\wedge \sum_{m\neq j}
(X_1,\dots,X_p)_n^n\wedge (Y_1,\dots,Y_q)^j_m +\\ + (-1)^{n+m}[X_n, Y_m]^c\wedge
(X_1,\dots,X_p)_n^n\wedge (Y_1,\dots,Y_p)_m^m  +\\ + \sum_{m,i}(-1)^{i+m}
[X_i,Y_m]^v\wedge \sum_{n\neq i} (X_1,\dots,X_p)^i_n\wedge (Y_1,\dots,Y_q)_m^m .
                                            \endmultline            $$
     On the other hand,
        $$ \multline
\dt[X,Y] = \sum_{i,j} (-1)^{i+j}[X_i,Y_j]^c\wedge (X_1,\dots,X_p)_i^i
\wedge (Y_1,\dots,Y_q)_j^j +\\
+ \sum_{i,j} (-1)^{i+j}[X_i,Y_j]^v\wedge\sum_{m\neq i} (X_1,\dots,X_p)_m^i
\wedge (Y_1,\dots,Y_q)_j^j +\\ + \sum_{i,j} (-1)^{i+j}[X_i,Y_j]^v\wedge \sum_{n\neq j}
(X_1,\dots,X_p)_i^i \wedge (Y_1,\dots,Y_q)_n^j
                                            \endmultline            $$
\noindent
and the required equality follows.
    \endproof
The rules of contractions for lifts of forms and vector fields we list in the
following Proposition.

    \proclaim{Proposition 2.6}{}
For $\zm\in \zF^1(\ezp_M)$ and $X\in \zF^1(\ezt_M)$ we have
        $$ \split \langle \vt(\zm),\vt(X)\rangle  &= \langle \ezt_M^{\textstyle *}
\zm, X^v\rangle  =0, \\
    \langle \dt\zm,\vt(X)\rangle = \langle \vt\zm,\dt(X)\rangle  &= \vt(\langle
\zm,X\rangle ) = \ezt_M^{\textstyle *} \langle \zm,X\rangle, \\
    \langle \dt\zm,\dt X\rangle &= \dt (\langle\zm, X\rangle ).
                                            \endsplit               $$
    \endproclaim
 These formulas can be easily generalized for $X\in \zF^r(\ezt_M)$ (or $\zm \in\zF^r
(\ezp_M)$). We obtain
        $$\split \langle \dt X,\vt\zm\rangle \overset \text(def) \to = \iv \dt X
&= \vt (\text{\tenrm i}_\zm X),\\
    \langle \dt X, \dt \zm \rangle \overset \text(def) \to  = \id \dt X &= \dt
(\text{\tenrm i}_\zm X),\\
    \langle \vt(X), \vt(\zm)\rangle \overset \text(def) \to = \iv \vt X&=0
                                \endsplit                   \tag 2.9    $$
for $X\in \zF^r(\ezt_M),\ \zm\in \zF^1(\ezp_M)$
and also
        $$\langle \dt\zm , \dt X\rangle =0 \ \ \text{for}\ \  \zm\in
\zF^r(\ezp_M),\ \ X\in\zF^r(\ezt_M), r>2.
                                                                $$

\subhead 3. Lie derivations of forms and multivector fields \endsubhead

The derivation $\dt$ on forms is strictly related to the Lie derivation in the algebra of
differential forms. It can be considered as a total Lie
derivative. This point of view is justified by Proposition~3.1.
Theorem~3.2 gives an analogous result for $\dt$ on multivector
fields.
    \proclaim{Proposition 3.1}{\cite{PiTu}} Let $X\colon M\rightarrow \sT M$ be a vector
field on $M$. Then we have
        $$ \ll_X\zm = X^\*\dt\zm .
                                                                \tag 3.1$$
        \endproclaim
    \proof Since $X^\*\colon \zF(\ezp_{\sss TM})\rightarrow \zF(\ezp_M)$ is a
homomorphism and $\dt \colon \zF(\ezp_M)\rightarrow \zF(\zp_{\sss TM})$ is a
$\vt$-derivation, the mapping
        $$ X^\*\dt\colon \zF(\ezp_M)\rightarrow \zF(\ezp_M)
                                                                \tag 3.2 $$
     is also a derivation.
It follows that it is enough  to verify the formula for functions and their
differentials. We have for a function $f$ on $M$
    $$ X^\*\dt f = \dt f \circ X = \langle X, \xd f\rangle  = X(f)   $$
and
    $$  X^\*\dt \xd f = X^\*\xd \dt f = \xd (X^\*\dt f) = \xd(X(f)) = \ll_X \xd f .  $$
    \endproof

    To get an analogue to this proposition in the case of vector
     fields, we need  an
analogue of the pull-back of forms with respect to a vector field.
Let $X$ be a vector field on $ M$. We have the decomposition of $\sT_{X(M)}\sT M$ into
horizontal (tangent to $X(M)$ )    and vertical parts
        $$ \sT_{X(M)}\sT M = H_X\sT M\oplus V_X\sT M .
                                                                \tag 3.3 $$

The canonical projection  $\sT_X\sT M \rightarrow V_X\sT M$ we denote by $P_X$. Let
$Y\in \zF^1(\ezt_{\sss T M})$ be a vector field on $\sT M$. By $X^\* Y$ we denote the
unique vector field on $M$ such that
        $$ \vt(X^\*Y)|_{X(M)} = P_XY .
                                                                \tag 3.4 $$
     The mapping $Y\mapsto  X^\*Y$ has the unique extension to a morphism of graded
algebras
        $$  X^\* \colon \zF(\ezt_{\sss TM}) \rightarrow \zF(\ezt_M) .
                                                                \tag 3.5 $$

    \proclaim{Theorem 3.2}{}    Let $X\colon M\rightarrow \sT M$ be a vector
field on $M$. For each $\zl\in \zF(\ezt_M)$ we have the following formulae:
        $$ \ll_X\zl = X^\*\dt\zl $$
and
        $$ \zl = X^{\textstyle *} \vt\zl.
                                                                $$
    \endproclaim
    \proof Since $X^\*\colon \zF(\ezt_{\sss TM})\rightarrow \zF(\ezt_M)$ is a
homomorphism and $\dt \colon \zF(\ezt_M)\rightarrow \zF(\ezt_{\sss TM})$ is a
$\vt$-derivation the mapping
        $$ X^\*\dt\colon \zF(\ezt_M)\rightarrow \zF(\ezt_M)
                                                                \tag 3.6$$
     is also a derivation.
It follows that it is enough  to verify the formula for functions and vector fields.
We have for a function $f$ on $M$ and a vector field $Y$ on $M$
    $$ X^\*\dt f = \dt f \circ X = \langle X, \xd f\rangle  = X(f)=\ll_Xf   $$
and (Proposition 2.6)
        $$ \ezt^{\textstyle *} _M(X^{\textstyle *} \dt Y(f)) = \vt((X^{\textstyle
*} \dt Y)(\dt f)) = \vt(X^{\textstyle *} \dt Y)(\dt f).
                                                                $$
On the other hand,
        $$ (P_X\dt Y)\circ X = \dt Y \circ X - \sT X\circ Y
                                                                $$
and, consequently,
        $$ (X^{\textstyle *} \dt Y)(f) = \langle \dt Y, \dt \xd f\rangle \circ X
-\langle  \sT X\circ Y, \dt\xd f\rangle .
                                                                $$
Since
        $$ \langle \dt Y, \dt \xd f\rangle \circ X = (\dt(\langle Y,\xd f\rangle ))\circ X
= X(Y(f))
                                                                $$
and
        $$ \langle  \sT X\circ Y, \dt\xd f\rangle = \langle Y, X^{\textstyle *}
\xd\dt f\rangle = \langle Y,\xd X(f)\rangle = Y(X(f)),
                                                                $$
we get
        $$ (X^{\textstyle *} \dt Y)(f) = X(Y(f)) - Y(X(f)) = [X,Y](f)= (\ll_XY)(f).
                                                                $$
The second identity follows directly from the definition.
    \endproof

    \specialhead Example. \endspecialhead

Let be $M=\R$, $X = X(x)\dfrac{\partial }{\partial x} $ and $Y =
Y(x)\dfrac{\partial }{\partial x}$. We have then
    $$\dt Y = Y(x)\dfrac{\partial }{\partial x} +
Y'(x)\dot x \dfrac{\partial }{\partial\dot x} $$
    and
    $$\dt Y|_{X} = Y(x)\dfrac{\partial }{\partial x} +
Y'(x)X(x) \dfrac{\partial }{\partial\dot x}. $$
The subspace of horizontal vectors is spanned by
$\dfrac{\partial }{\partial x} + X'(x) \dfrac{\partial
}{\partial \dot x} $. Hence the decomposition of $\dt Y|_X$ into
horizontal and vertical parts is the following:
        $$\multline Y(x)\dfrac{\partial }{\partial x} +
Y'(x)X(x) \dfrac{\partial }{\partial\dot x} =\\
    = Y(x)\left( \frac{\partial }{\partial x} +
X'(x)\frac{\partial }{\partial \dot x} \right) + (X(x)Y'(x)
-Y(x)X'(x))\frac{\partial }{\partial \dot x}.\endmultline $$
    The vertical part is obviously the vertical lift of
$$(X(x)Y'(x)-Y(x)X'(x))\dfrac{\partial }{\partial x} = [X,Y](x)=
\ll_X Y (x).$$
\smallskip

    In the following sections we represent a Poisson structure
by a vector bundle morphism $\sT^{\textstyle *} M\rightarrow \sT
M$ rather than by a bi-vector field. We need then a formula for
the Lie derivative of a bi-vector field $\zl$ expressed in terms
of $\widetilde{\zl}^1$. In order to get it, we  first identify an
operation dual to $X^{\textstyle *} $.

    The  projection $ \sT\ezt_M$ induces an isomorphism $\sT\ezt_M \colon H_X \sT M
\rightarrow \sT M$ and the dual isomorphism $(\sT\ezt_M)'\colon \sT^{\textstyle *}
M\rightarrow H_X^{\textstyle *} \sT M $.

When composed with $(I-P_X)'$, this isomorphism gives an isomorphism
$\sT^{\textstyle *} M\rightarrow V^\circ _X\sT M$, where $V^\circ _X\sT M$ is the annihilator
of $V_X\sT_M$ in $\sT^{\textstyle *} _{X(M)}\sT M$.
It follows that for $\zn\in \zF^1(\ezp_{\ssT M})$ the decomposition
        $$  \zn =  (\zn - \vt (X^{\textstyle *} \zn) ) + \vt(X^{\textstyle *} \zn)
                                                                $$
is, on $X(M)$, the decomposition related to the decomposition
        $$   \sT^{\textstyle *} _{X(M)}\sT M = H^\circ _X\sT
M\oplus V^\circ _X\sT M,
                                                                $$
     dual to (3.3).
For $\zn = \dt \zm$ we have $\vt(X^{\textstyle *} \zn )= \vt(\ll_X\zm)$
(Proposition 3.1). Hence, for $Y\in \zF^1(\ezt_{\ssT M})$ we have
        $$ \langle \zm, X^{\textstyle *} Y\rangle = (\langle \dt\zm,
\vt(X^{\textstyle *} Y)\rangle) \circ X = (\langle \dt\zm - \vt(\ll_X\zm), Y\rangle
)\circ X
                                                \tag 3.7            $$

It is convenient to introduce an operation

        $$X_{\textstyle +}\colon  \zF(\ezp_M)\rightarrow \zF(\ezp_{\ssT M}) \colon
\zm \longmapsto \dt\zm - \vt(\ll_X\zm),
                                                                     $$
which is a $\vt$ -derivation of degree $0.$


It follows directly from (3.7) that
for $Y\in \zF^r(\ezt_{\ssT M})$ and $\zm_1,\dots ,\zm_r\in
\zF^1(\ezp_M)$, $r>0$,  we have
        $$
          \langle (X_{\textstyle +}\zm_1)\wedge \dots\wedge (X_{\textstyle +}\zm_r),Y\rangle \circ X
= \langle \zm_1\wedge \dots\wedge\zm_r,
X^{\textstyle *} Y\rangle ,
                                                \tag 3.8            $$
so the mapping
$$X_{\textstyle *}(\zm_1\wedge\dots\zm_r)=(X_{\textstyle +}\zm_1)\wedge\dots\wedge (X_{\textstyle +}\zm_r)$$
regarded as a homomorhism of graded algebras
$$X_{\textstyle *}:\zF(\ezp_M)\rightarrow\zF(\ezp_{X(M)}),$$
where
$$\ezp_{X(M)}:\sT^\*_{X(M)}\sT M\rightarrow X(M),$$
is dual to $X^{\textstyle *}:$
$$\langle X_{\textstyle *}\zm ,Y\rangle\circ X
=\langle\zm ,X^{\textstyle *}Y\rangle .$$
For 1-forms we shall write $X_{\textstyle *}\zm$ instead of
$X_{\textstyle +}\zm.$
    \proclaim{Proposition 3.3}{} Let $X$ be a vector field on $M$ and let $\zl$ be a
2-vector field on $M$. $\ll_X \zl = 0$ if and only if for every pair $f,g $ of
functions on $M$
    $$ \langle  X_{\textstyle *} \xd g , \widetilde{\dt\zl}^1(X_{\textstyle *} \xd
f)\rangle =0. $$
    \endproclaim
    \proof
 From Theorem 3.2 it follows that $\ll_X\zl =0$ iff
        $$\langle \xd g, \widetilde{X^{\textstyle *}\dt\zl }^1
          (\xd f)\rangle =0
                                                                $$
    for each pair $(f,g)$ of functions on $M$.
  From (3.8) we have then that  $\ll_X\zl =0$ if and only if
        $$ \langle X_{\textstyle *}\xd g, \widetilde{\dt\zl}^1X_{\textstyle *} \xd
f \rangle = 0.
                                                                $$
%
%
    \endproof

\subhead 4. Symplectic and Poisson structures \endsubhead
    In this section we give definitions of symplectic and
Poisson structures represented by  morphisms of tangent and
cotangent bundles. These definitions do not make use of the  exterior
derivative and of the Schouten bracket.

By the standard definition, a symplectic structure on $M$ is a
nondegenerate, closed two-form $\zw$ on $M$. The canonical example is
the 2-form $\zW_M$ on $\sT^{\textstyle *} M$.
On the other hand, any 2-form $\zm$ on $M$ can be represented by
        $$ \widetilde{\zm}^i \colon   {\bigwedge}^i\sT M \rightarrow
{\bigwedge}^{2-i}\sT^{\textstyle *} M, \ \ i=0,1,2.
                                                                $$
The standard definition is expressed, in fact,  in terms of $\widetilde{\zm}^0$.
We can also formulate a definition of a symplectic structure in terms of
$\widetilde{\zm}^2$ making use of the well known formula for the exterior derivative:

    \proclaim{Proposition 4.1}{}
A linear fibre bundle morphism $ \zn \colon {\bigwedge}^2\sT M \rightarrow
{\bigwedge}^0\sT^{\textstyle *} M$ represents a symplectic structure on $M$ if it is
 nondegenerated (in the obvious meaning)  and
        $$  X(\zn(Y,Z)) + Y(\zn(Z,X)) + Z(\zn(X,Y)) - \zn([X,Y],Z) - \zn([Y,Z],X)
-\zn([Z,X],Y) =0
                                                                $$
      where $X,Y,Z \in \zF^1(\ezt_M)$.
    \endproclaim
     Now, we provide a definition of a (pre-)symplectic structure in terms of
$\widetilde{\zm}^1$.

First, we  need a condition for a 2-form $\zm$ to be closed  in terms of $\widetilde{\zm}^1$.

    \proclaim{Theorem 4.2}{}

A 2-form $\zm$ on $M$ is closed if and only if
        $$ \dt \zm = (\widetilde{\zm}^1)^{\textstyle *} \zW_M $$
    \endproclaim
    \proof
 From the formula 2.2 we have
        $$\xd\xit\zm = (\widetilde{\zm}^1)^{\textstyle *} \zW_M
                                                                $$
and, consequently,
        $$ \dt\zm = \xd\xit\zm + \xit\xd\zm = (\widetilde{\zm}^1)^{\textstyle *}
\zW_M - \xit\xd\zm.
                                                                $$
Since $\xit \xd\zm = 0$ iff $\xd\zm =0 $, we get that $\xd\zm = 0$ if and only if
        $$ \dt \zm = (\widetilde{\zm}^1)^{\textstyle *} \zW_M .
                                                                $$
    \endproof

A condition for a vector  bundle morphism $\zn \colon \sT M\rightarrow
\sT^{\textstyle *}M $ to define a symplectic structure on $M$ is given by the
following theorem.

    \proclaim{Theorem 4.3}{} An  isomorphism $\zn \colon \sT M \rightarrow
\sT^{\textstyle *}M $ of vector bundles defines a symplectic structure on $M$
if and only if the following diagram is commutative
    $$   \CD
            \sT\sT^\*M     @>\dsize \widetilde{\zW_M}^1 >>   \sT^\*\sT^\*M \\
            @A  \dsize \sT \zn AA   @V\dsize \sT^\*\zn VV\\
                \sT\sT M    @.  \sT^\*\sT M  \\
            @V \dsize \zk_M VV @A \dsize \za_M AA \\
              \sT\sT M  @>\dsize \sT\zn>> \sT\sT^\* M
                                                \endCD \ \ .
                                                    \tag    4.1     $$
    \endproclaim
    \proof
Commutativity of the diagram 4.1 is equivalent to the equality
        $$ \za_M\circ \sT \zn \circ \zk_M = \sT^{\textstyle *} \zn \circ
\widetilde{\zW_M}^1 \circ \sT \zn
                                                    \tag 4.2        $$
and to the dual (with respect to proper pairings) equality
        $$ \za_M\circ \sT \zn^{\textstyle *}  \circ \zk_M = \sT^{\textstyle *} \zn \circ
(\widetilde{\zW_M}^1)^{\textstyle *}  \circ \sT \zn.
                                                    \tag 4.3        $$
     Since $ (\widetilde{\zW_M}^1)^{\textstyle *} = - \widetilde{\zW_M}$, we get from
4.3 that
        $$ \sT\zn = - \sT \zn^{\textstyle *}
                                                                $$
and, consequently, that $\zn^{\textstyle *} = -\zn$. It proves that $\zn$ is
skew-symmetric and, consequently, that the exist a 2-form $\zm$ on $M$, such that
$\zn  = \widetilde{\zm}^1$.
Equality 4.2 reads now
        $$ \xd \zm = (\widetilde{\zm}^1)^{\textstyle *} \zW_M.
                                                                $$
It follows from Theorem 4.2 that $\zm$ is closed.
        \endproof

As in the case of symplectic structures, a Poisson structure on $M$ can be
represented by one of three equivalent objects:
\roster
    \item a bivector field $\zL$,
    \item a homomorphism $\widetilde{\zL}^1 $ of vector bundles,
    \item a bilinear, skew-symmetric function $\widetilde{\zL}^2$.
\endroster

The condition for $\zL$ to be a Poisson structure is vanishing of the Schouten
bracket: $[\zL,\zL]= 0$.

A skew-symmetric function
        $$\zl \colon {\bigwedge}^2 \sT^{\textstyle *} M \rightarrow
{\bigwedge}^0\sT M \simeq M\times \R
                                                                $$
defines a Poisson structure if the Jacobi identity is fulfilled:
        $$ \zl(\xd f,\zl(\xd g,\xd h)) + \zl(\xd g,\zl(\xd h,\xd f)) +
\zl(\xd h,\zl(\xd f,\xd g)) =0,
                                                                $$
where $f,g, h \in \zF^0(\ezp_M)$.

In order to get a  definition of a  Poisson structure in terms of a vector bundle
morphism $\zn \colon \sT^{\textstyle *} M\rightarrow \sT M$, we need the following
theorem.

    \proclaim{Theorem 4.4}{}

A 2-vector field $\zL$ on $M$ defines Poisson structure if and
only if
        $$ (\widetilde{\zL}^1)_{\textstyle *}\zL_M \subset \dt \zL  $$
(one can say that $\zL_M$ and $\dt\zL$ are $\widetilde{\zL}^1$-related),
where $\zL_M = (\zW_M)^{-1} $ is the canonical 2-vector field on $\sT^{\textstyle *}
M$.
    \endproclaim
    \proof
    In local coordinates $\zL = \dfrac{1}{2} \zL^{ij}\dfrac{\partial}{\partial x^i
}\wedge \dfrac{\partial}{\partial x^j } $.
The 2-vector field $\zL$ defines a Poisson structure if and
only if
        $$ \sum_i(\zL^{ij}\dfrac{\partial}{\partial x^i }\zL^{kl} +
\zL^{ik}\dfrac{\partial}{\partial x^i } \zL^{lj} + \zL^{il}\dfrac{\partial}{\partial
x^i } \zL^{jk}) =0 .
                                                                \tag 4.4$$
     The canonical 2-vector field $\zL_M$ on $\sT^\*M$ is given by
    $$ \zL_M = \sum_i\dfrac{\partial}{\partial x^i }\wedge \dfrac{\partial}{\partial p_i }.$$
    Moreover,
    $$ \dt\zL = \frac{1}{2}
     \sum_{ijk}\frac{\partial\zL^{ij}}{\partial x^k}\dot
x^k\dfrac{\partial}{\partial\dot x^i } \wedge \dfrac{\partial}{\partial\dot x^j } +
\sum_{ij} \zL^{ij}\dfrac{\partial}{\partial\dot x^i }\wedge \dfrac{\partial}{\partial x^j }.
                                                                \tag 4.5 $$
    For $a\in \sT^\*M,$
    $$ \dot x^i((\widetilde{\zL}^1)(a))= \sum_j\zL^{ji}p_j(a)$$
and, consequently,
    $$ \aligned (\widetilde{\zL}^1)_\*(\dfrac{\partial}{\partial x^i }) &=
\dfrac{\partial}{\partial x^i } + \sum_{j,k}p_j(\dfrac{\partial}{\partial x^i }\zL^{jk})
\dfrac{\partial}{\partial\dot x^k } \\
        (\widetilde{\zL}^1)_\*(\dfrac{\partial}{\partial p_i }) &=
\sum_j\zL^{ij}\dfrac{\partial}{\partial\dot x^j }.\endaligned $$

Hence
        $$(\widetilde{\zL}^1)_\*\zL_M = \sum_{i,j}\zL^{ij}\dfrac{\partial}{\partial x^i
}\wedge  \dfrac{\partial}{\partial\dot x^j } +
\sum_{ijk}p_l\zL^{ij}(\dfrac{\partial}{\partial x^i }\zL^{lk})
\dfrac{\partial}{\partial\dot x^k }
\wedge \dfrac{\partial}{\partial\dot x^j } .
                                                                \tag 4.6$$
    It follows from 4.5 and 4.6 that $ (\widetilde{\zL}^1)_{\textstyle
*}\zL_M \subset \dt \zL  $ if and only if
    $$ \zL^{lk}p_l\dfrac{\partial}{\partial x^k }\zL^{ij} =
\zL^{kj}(\dfrac{\partial}{\partial x^k }\zL^{li})p_l -
\zL^{ki}(\dfrac{\partial}{\partial x^k }\zL^{lj})p_l, $$
    but this is equivalent to 4.4.
    \endproof

Now, we are ready for a proof of an analogue of Theorem 4.3.

    \proclaim{Theorem 4.5}{} A vector bundles morphism  $\zl \colon \sT^{\textstyle
*} M \rightarrow \sT M$ defines a Poisson structure on $M$ if and only if the
following diagram is commutative

    $$   \CD
            \sT\sT^\*M     @<<\dsize \zL_M<   \sT^\*\sT^\*M \\
            @V  \dsize \sT\zl VV    @A\dsize \sT^\*\zl AA\\
                \sT\sT M    @.  \sT^\*_{im(\zl)}\sT M  \\
            @A \dsize \ze_M' AA @V \dsize \ze_M VV \\
              \sT\sT M  @<<\dsize \sT\zl< \sT\sT^\* M
                                                \endCD\ \  .
                                                    \tag    4.7     $$

    \endproclaim
\proof
Commutativity of the diagram is equivalent to
        $$ \sT\zl \circ \widetilde{\zL_M}^1 \circ  \sT^{\textstyle *} \zl = \ze_M'\circ \sT\zl
\circ \ze_M
                                                                $$
and to the dual equality
        $$ \sT\zl \circ (\widetilde{\zL_M}^1)^{\textstyle *}  \circ  \sT^{\textstyle *} \zl =
\ze_M'\circ \sT\zl^{\textstyle *}  \circ \ze_M.
                                                                $$
Since $(\widetilde{\zL_M}^1) ^{\textstyle *} = -\widetilde{\zL_M}^1$ we get
$\sT\zl^{\textstyle *} = -\sT\zl$ and, consequently, $\zl^{\textstyle *}  = -\zl $.
It follows that there exists a bivector field $\zL$ on $M$ such that  $\zl =
\widetilde{\zL}^1$. From Theorem 4.4 we have that $\zL$ is a Poisson structure if and
only if the diagram 4.7 is commutative.
\endproof

In the diagrams, $\sT^\* \zl$
is a relation (not a mapping) and diagrams
are in the category of relations.

In order to illustrate  the condition 4.7, let us consider the case of a linear Poisson
structure. Let $M= V$ be a vector space. We have obvious identifications:
        $$ \split \sT V =V\times \sT_0V =V\times V, & \quad \sT^{\textstyle *} V =
V\times V^{\textstyle *}, \\
        \sT\sT V = (V\times V)\times (V\times V), &\quad \sT\sT^{\textstyle *} V = (V\times
V^{\textstyle *} )\times (V\times V^{\textstyle *} ), \\
    \sT^{\textstyle *} \sT^{\textstyle *} V = (V\times V^{\textstyle *} )\times
(V^{\textstyle *} \times V) .
                                \endsplit                           $$
With these identifications the canonical morphisms $\ze_M,\zk_M,
\widetilde{\zL_M}^1$
look like follows:
        $$ \split
\ze_V\colon \sT^{\textstyle *} \sT V\ni (v,w;a,b) &\longmapsto (v,b;w,a) \in
\sT\sT^{\textstyle *} V, \\
    \zk_V\colon \sT\sT V\ni (v,w;x,y)&\longmapsto (v,x;w,y)\in \sT\sT V, \\
    \widetilde{\zL_V}^1\colon \sT^{\textstyle *} \sT^{\textstyle *} V \ni
(v,a;w,b)& \longmapsto (v,a; -w,b)\in \sT\sT^{\textstyle *} V.
                                \endsplit                       $$
A linear Poisson structure $\zL$ is given by a mapping
        $$ \widetilde{\zL}^1\colon \sT^{\textstyle *} V\ni (v,a) \longmapsto
(v,\zl(v,a))\in \sT V,
                                                                $$
where $\zl\colon V\times  V^{\textstyle *} \rightarrow V $ is bilinear.
We have for such $\zL$:
        $$\split \sT\widetilde{\zL}^1\colon (v,a;w,b)&\longmapsto (v,\zl(v,a); w,
\zl(v,b) + \zl(w,a)),\\
    \sT^{\textstyle *} \widetilde{\zL}^1 \colon (v,\zl(v,d); a,b) &\longmapsto
(v,d; a + {}^{\textstyle *} \zl(b,d), \zl^{\textstyle *} (v,b) ),
                                \endsplit                           $$
where ${}^{\textstyle *} \zl$ and $\zl^{\textstyle *} $ are conjugate to $\zl$ with
respect to the left and to the right argument respectively.
The condition 4.7 reads as follows:

The mapping $\widetilde{\zL}^1$ defines a Poisson structure if and only if the following conditions are
satisfied for each $v,a,b$:
        \roster
    \item $ \zl(v,b) = - \zl^{\textstyle *} (v,b)$,
    \item $\zl(v,{}^{\textstyle *}\zl(b,a) ) = \zl(\zl(v,a),b) -
     \zl(\zl(v,b),a),$
        \endroster
which means that
$${}^{\textstyle *}\zl\colon V^{\textstyle *}\times V^{\textstyle *}\rightarrow V^{\textstyle *}$$
is a Lie bracket.

 \subhead 5. Tangent  Poisson structures  \endsubhead

Let $(M,\zL)$ be a Poisson manifold, where $M$ is a manifold and $\zL\in
\zF^2(\ezt_M)$ satisfies $[\zL,\zL] = 0$. It follows from Theorem 2.5 that $(\sT M,\dt \zL)$
is also a Poisson manifold:
        $$ [\dt\zL, \dt\zL] = \dt[\zL,\zL] =0.
                                                                $$
     We call $\dt\zL$ the {\it tangent Poisson structure}.
Let be, in local coordinates,
        $$ \zL = \frac{1}{2} \zL^{ij}(x)\frac{\partial }{\partial x^{i}}\wedge
\frac{\partial }{\partial x^{j}}, \ \ \  \zL^{ij} = -\zL^{ji},
                                                                $$
then
        $$\dt\zL = \zL^{ij}(x) \frac{\partial }{\partial  x^{i}} \wedge
\frac{\partial }{\partial \dot x^{j}} + \frac{1}{2}\frac{\partial
\zL^{ij} }{\partial x^{k}} \dot x^k \frac{\partial }{\partial\dot
x^{i}}\wedge \frac{\partial }{\partial \dot x^{j}} .
                                        \tag 5.1                    $$
The Poisson structure $\zL$ defines a Poisson bracket $\{f,g\}_{\zL} = \zL(\xd f\wedge \xd
g)$ which provides the algebra $\zF^0(\ezt_M)$ of smooth functions with a structure of
Lie algebra. The tangent Poisson structure defines a Poisson bracket $\{,\}_{\sdt\zL}$ on $\sT M$ which
is characterized by the following relations:
        $$ \split \{\dt f, \dt g \}_{\sdt\zL} &= \dt\{f,g\}_\zL ,\\
        \{\dt f, \vt g \}_{\sdt\zL}= \{\vt f, \dt g\} &= \vt\{f,g\}_\zL , \\
        \{\vt f, \vt g \}_{\sdt\zL} &= 0 .
                                        \endsplit \tag 5.2          $$
    This is exactly the lift of $\zL$ defined in \cite{SdA,Co3}.
For a 1-form $\zm\in \zF^1(\ezp_M)$ we put $\zL_\zm = \xi_\zm \zL$. For a function
$f$ the  vector field $\zL_{\hbox{\sevenrm d}f}$ is the {\it hamiltonian vector field
generated by $f$}.
    \proclaim{Theorem 5.1}{}
The tangent Poisson structure $\dt\zL$ is linear with respect to the vector bundle
structure in $\ezt_M \colon \sT M \rightarrow M$. Moreover, for $\zm,\zy \in
\zF^1(\ezp)$ we have the following formula:
        $$ \{\xit\zm,\xit \zy \}_{\sdt \zL} = \xit (\xd\langle \zL,\zm\wedge
\zy\rangle  + \xi_{\zL_\zm}\xd\zy - \xi_{\zL_\zy}\xd \zm ).
                                                \tag 5.3            $$
    \endproclaim
    \proof
It enough to proof the formula 5.3. We have
        $$ \multline
    \{ \xit \zm, \xit \zy\}_{\sdt\zL} = \langle \dt\zL, \xd \xit\zm\wedge \xd \xit
\zy \rangle = \\
    = \langle \dt\zL , \dt \zm\wedge \dt\zy - \dt\zm \wedge  \xit \xd \zy - \xit \xd
\zm \wedge \dt\zy + \xit\xd \zm \wedge \xit\xd \zy \rangle . \endmultline
                                                \tag 5.4            $$
    Since the 1-forms $\xit\xd\zm$ and $\xit\xd\zy$ are vertical, we have
        $$ \langle \dt \zL, \xit \xd \zm \wedge \xit \xd \zy \rangle =0.
                                                                $$
     Moreover,
        $$
    \langle \dt \zL,  \dt \zm\wedge \dt \zy \rangle  = \xi_{\sdt
\zy}\xi_{\sdt \zm}\dt\zL
    = \xi_{\sdt\zy}\dt(\xi_\zm\zL) = \dt(\langle \zL,\zm\wedge \zy \rangle )
                                                                $$
and
        $$ \multline
    \langle \dt \zL,  \dt \zm\wedge \xit\xd \zy \rangle  =\langle \dt(\zL_\zm),\xit
\xd\zy \rangle = \xi_{\sdt(\zL_ \zm)}\xi_\zG\vt (\xd\zy)  =\\
    = -\xi_\zG\xi_{\sdt\zL_\zm}\vt \xd\zy = -\xi_\zG\vt(\xi_{\zL_\zm}\xd\zy ) =
-\xit (\xi _{\zL_\zm}\xd\zy),
                                        \endmultline
                                                                $$
where $\zG$ is a second order vector field (see Section 2).
Similarly,
        $$ \langle \dt\zL, \xit \xd \zm \wedge \dt \zy \rangle  =
\xit(\xi_{\zL_\zy}\zm)
                                                                $$
and the theorem follows.
    \endproof

Since $\dt \zL$ is a linear Poisson structure, it defines an algebroid structure in
the dual vector bundle $\sT^{\textstyle *} M$. The theorem provides an explicit
formula for the Lie bracket in $\zF^1(\ezp_M)$:
        $$ \split \{\zm,\zy\}_\zL &= \xd\langle \zL, \zm\wedge \zy\rangle  +
\xi_{\zL_\zm}\xd\zy - \xi_{\zL_\zy}\xd\zm \\
                &= \ll_{\zL_\zm}\zy - \ll_{\zL_\zy}\zm - \xd\langle \zL,
\zm\wedge \zy \rangle .
                                        \endsplit                   $$

\noindent
This is exactly the well-known extension of a Poisson bracket to
1-forms. Thus we have got a new way to define it using the tangent lift.

\specialhead Example.\endspecialhead

Let us consider the Poisson structure on $\R^2$ given by

    $$\zL = x^2\frac{\partial}{\partial x}\wedge\frac{\partial}{\partial
y}.$$

The tangent lift of this structure is given by the following formula:

$$\dt\zL  =x^2(\frac{\partial}{\partial \dot x}\wedge
\frac{\partial}{\partial y} +
\frac{\partial}{\partial x}\wedge\frac{\partial}{\partial \dot y})  +
2x\dot x \frac{\partial}{\partial \dot x}\wedge
\frac{\partial}{\partial \dot y}.$$

For   1-forms   $\zm    =\zm_1(x,y)\xd x+\zm_2(x,y)\xd y$    and    $\zy
=\zy_1(x,y)\xd x+\zy_2(x,y)\xd y$ we have  $\xit\zm  =\zm_1(x,y)\dot  x+
\zm_2(x,y)\dot y$ and  $\xit\zy  =\zy_1(x,y)\dot  x+\zy_2(x,y)\dot
y$. We can easily calculate the Poisson bracket

$$\multline \{\xit\zm ,\xit\zy\}_{\dt\zL}=
x^2\left[ \zm_1\left(
\frac{\partial  \zy_1}{\partial   y}\dot   x+\frac{\partial   \zy_2}
{\partial y}\dot y\right) -
\zm_2\left(
\frac{\partial  \zy_1}{\partial  x}\dot   x+\frac{\partial   \zy_2}
{\partial x}\dot y\right)\right. \\
    \left. - \zy_1\left(
\frac{\partial  \zm_1}{\partial   y}\dot   x+\frac{\partial   \zm_2}
{\partial y}\dot y\right) +
\zy_2\left(
\frac{\partial  \zm_1}{\partial   x}\dot   x+\frac{\partial   \zm_2}
{\partial x}\dot y\right)\right] +
2x\dot x(\zm_1\zy_2-\zm_2\zy_1).\endmultline$$

Hence,

$\{\xit\zm ,\xit\zy\}_{\dt\zL}=
\xit\{\zm ,\zy\}_\zL,$ where

$$\multline \{\zm ,\zy\}_\zL =\left[  x^2\left(  \zm_1\frac{\partial  \zy_1}
{\partial y}-
\zm_2\frac{\partial  \zy_1}
{\partial x}-
\zy_1\frac{\partial  \zm_1}
{\partial y}+
\zy_2\frac{\partial  \zm_1}
{\partial x}\right) +
2x(\zm_1\zy_2-\zm_2\zy_1)\right] \xd x \\+
x^2\left(
\zm_1\frac{\partial  \zy_2}
{\partial y}-
\zm_2\frac{\partial  \zy_2}
{\partial x}-
\zy_1\frac{\partial  \zm_2}
{\partial y}+
\zy_2\frac{\partial  \zm_2}
{\partial x}\right) \xd y.\endmultline$$

In   particular,   $\{   x\xd y,\xd y\}_\zL    =    x^2\xd y$    and    $\{
x\xd y,\xd x\}_\zL=-2x^2\xd x$.

\subhead 6. Canonical vector fields \endsubhead

It is well known  \cite{Tu2}  that for a symplectic manifold $(M,\zw)$ the tangent structure
$(\sT M,\dt\zw)$ is also a symplectic manifold.
We use Proposition 3.1,   to get a simple proof that canonical vector fields on a
symplectic manifold are Lagrangian submanifolds with respect to the tangent
symplectic structure (compare with \cite{SdA}). It justifies the
concept of a generalized canonical system as a Lagrangian
submanifold of the tangent Poisson manifold.
        \proclaim{Proposition 6.1}{}  Let $X\colon M\rightarrow \sT M$ be a vector field
on a symplectic manifold $(M,\zw)$. The vector field $X$ is canonical, i.~e., $\ll_X\zw =0$ if and
only if $X(M)$ is a Lagrangian submanifold of $(\sT M, \dt \zw)$.

    \endproclaim
    \proof From the Proposition ~3.1,  $\ll_X\zw = X^\*\dt \zw.$ Hence
$\ll_X\zw=0$ if and only if $X^\*\dt\zw=0$, but the last means that $X(M)$ is
isotropic in  $(\sT M, \dt \zw)$ and, consequently, Lagrangian (the dimension of
$X(M) $ is equal to the dimension of $M$).

    \endproof

    We have an analogue of this theorem for Poisson manifolds. To formulate it we
need the definition of a Lagrangian submanifold of a Poisson manifold.
    \proclaim{Definition}{} Let $(M,\zL) $  be a Poisson manifold. A submanifold
$N\subset M$ is {\rm Lagrangian } if for each $m\in N$
        $$ \widetilde{\zL}^1(\sT_m N)^o =\widetilde{\zL}^1(\sT^\*_mM)\cap \sT_mN
                                                         \tag 6.1 $$

    \endproclaim

        \proclaim{Theorem 6.2}{}  Let $X\colon M\rightarrow \sT M$ be a vector field
on a Poisson manifold $(M,\zL)$. $X$ is canonical, i.~e., $\ll_X\zL =0$ if and
only if $X(M)$ is a Lagrangian submanifold of $(\sT M, \dt \zL)$.

    \endproclaim
    \proof

    From Proposition 3.3, $\ll_X\zL=0$ if and only if, for $f,g\in C^\infty(M)$,
    $$ \langle \widetilde{\dt\zL}^1(X_{\textstyle *} \xd f), X_{\textstyle *} \xd g
\rangle =0.
                                                                $$

    On the other hand, $X(M)$ is Lagrangian if for each $v\in X(M)$

        $$  \{ \widetilde{\dt\zL}^1(X_{\textstyle *} \xd f)(v)| f\in C^\infty(M) \} = \sT_vX(M)
\cap \widetilde{\dt\zL}^1 (\sT^\*_v\sT M).
                                                        \tag 6.2     $$
    Since $w\in \sT_v X(M)$ is characterized by $X_{\textstyle *} \xd f(w)=0 $, the
inclusion $''\subset '' $   is equivalent to
     $$ \langle \widetilde{\dt\zL}^1(X_{\textstyle *} \xd f)(v), X_{\textstyle *}
\xd g \rangle =0.
                                                                $$
Thus, if $X(M)$ is Lagrangian then $X$ is canonical.
In order to prove that for a canonical $X$ the submanifold  $X(M)$ is
Lagrangian, it is
enough to show that the inclusion $''\supset'' $ holds.
    Let be $a\in \sT^\* _v\sT M$. There are functions $f,g$ on $M $ such that $a =
\xd_v(\dt f) + \xd_v(\vt g)$. We have then (2.9)
    $$ \widetilde{\dt\zL}^1(a) = \dt(\widetilde{\zL}^1\xd f)
+\vt(\widetilde{\zL}^1\xd g).
                                                    \tag 6.3        $$
    It follows that
    $$\sT\ezt_M (\widetilde{\dt\zL}^1a) = \sT\ezt_M(\dt(\widetilde{\zL}^1\xd f) =
\widetilde{\zL}^1\xd_m f,
                                                        \tag 6.4    $$
where $m = \ezt_M v $.
    Thus the tangent projections of vectors in $\sT X(M)\cap \widetilde{\dt \zL}^1
(\sT^\* \sT M) $  are in the image of $\widetilde{\zL}^1$. For a vector
$w\in\widetilde{\zL}^1(\sT^\* M)$ there is only one  vector $v\in \sT X(M)$  such
that  $\sT\ezt_M v = w$, namely,  $v = \sT X(w)$. Let be  $w= \widetilde{\zL}^1\xd_mf$.
The inclusion $''\subset ''$ implies that     $ \widetilde{\dt\zL}^1(\xd_{X(m)} (\dt
f-\vt (Xf))$  is tangent to $X(M)$. Since (6.3 and 6.4)
        $$ \sT\ezt_M \widetilde{\dt\zL}^1(\xd (\dt f-\vt (Xf))(X(m))) =
\widetilde{\zL}^1\xd_mf =w,
                                                                $$
we have $\widetilde{\dt\zL}^1(X_{\textstyle *} \xd f)(X(m)) = \sT X(w)$.
    \endproof

\subhead 7. Tangent  Poisson-Lie structures \endsubhead

    The growing interest in Poisson-Lie structures justifies
analysis of every aspect of these structures. We show that the
tangent lift of a Poisson-Lie structure is a Poisson-Lie
structure (Theorem~7.1). Vertical and complete lifts of left
(right) invariant vector fields on $G$ turn out to be left
(right) invariant on $\sT G$.

    Let $M,\ N$ be differentiable manifolds. Since there is a
canonical identification $\sT(M\times N) \simeq \sT M\times \sT
N$, we have also canonical inclusions ${\bigwedge}^2 \sT M
\times {\bigwedge}^2 \sT N \subset {\bigwedge}^2 \sT (M\times
N)$ and $\zF^2(\ezt_M)\oplus \zF^2(\ezt_N) \subset
\zF^2(\ezt_{N\times M})$. It is trivial task to verify, that if
$ (M,\zL)$ and $(N,\zP)$ are Poisson manifolds then $(M\times N,
\zL \oplus \zP)$ is also a Poisson manifold. The Poisson
structure $\zL \oplus \zP$ is called the {\it product Poisson
structure}. Finally, we recall the notion of a {\it Poisson
map}. A smooth mapping $\zf \colon M \rightarrow N$ is a Poisson
map if
        $$ {\bigwedge}^2 \sT \zf\circ \zL = \zP\circ \zf
                                                                $$
     or, equivalently,
        $$\sT\zf \circ \widetilde{\zL} \circ \sT^{\textstyle *} \zf = \widetilde{\zP},
                                                                $$
where $\widetilde{\zL} = \widetilde{\zL}^1$, $\widetilde{\zP} = \widetilde{\zP}^1$.

Let $(G,m)$ be a Lie group and let $(G,\zL)$ be a Poisson manifold. We say that
$(G,m,\zL)$ is a {\it Poisson-Lie group} if the group multiplication $m$ is  a
Poisson mapping:
        $$ m\colon (G\times G, \zL \oplus \zL)\rightarrow (G,\zL)
                                                                $$
     or, equivalently,
        $$ \sT m \circ (\widetilde{\zL}\times \widetilde{\zL}) \circ
\sT^{\textstyle *} m = \widetilde{\zL}.
                                                    \tag 7.1        $$
We say in this case that $\zL$ is {\it multiplicative}.

Applying the tangent functor to 7.1, we get
        $$ \sT\sT m\circ ( \sT\widetilde{\zL}\times \sT\widetilde{\zL})\circ
\sT\sT^{\textstyle *} m = \sT \widetilde{\zL}.
                                                        \tag 7.2    $$
It is well known that $(\sT G, \sT m)$ is a Lie group.
    \proclaim{Theorem 7.1}{}
Let $(G,m,\zL)$ be a Poisson-Lie group. Then $(\sT G, \sT m, \dt \zL )$ is also a
Poisson-Lie group.
    \endproclaim
    \proof
We have to show that $\dt \zL$ is multiplicative with respect to $\sT m$, i.~e., that
        $$\sT\sT m\circ (\widetilde{\dt\zL}\times \widetilde{\dt\zL})\circ
\sT^{\textstyle *} \sT m = \widetilde{\dt\zL}.
                                                                $$
Since, due to (2.6),  $\widetilde{\dt \zL}= \zk_G\circ \sT\widetilde{\zL} \circ \ze_G$, it
follows from  functorial properties of $\zk$ and $\ze$ that
        $$ \split
\sT\sT m\circ (\widetilde{\dt\zL}\times \widetilde{\dt\zL}) \circ
\sT^{\textstyle *} \sT m &= \sT\sT m \circ  \zk_{G\times G}\circ (\sT \widetilde{\zL}
\times \sT \widetilde{\zL}) \circ \ze_{G\times G} \circ  \sT^{\textstyle *} \sT m \\
    &= \zk_G\circ \sT\sT m \circ (\sT \widetilde{\zL}\times \widetilde{\zL})\circ
\sT \sT^{\textstyle *} m \circ \ze_G
                                        \endsplit                   $$
and by (7.2) the required identity follows.
    \endproof

Let \fr g be the Lie algebra of the Lie group $G$. The tangent bundle $\sT G$ can be
trivialized by  right  or left translations:
        $$(K^r, \ezt_G)\colon \sT G\rightarrow \fr g\times G
                                                                $$
      or
        $$(\ezt_G ,K^l )\colon \sT G\rightarrow  G\times \fr g,
                                                                $$
where $K^r(K^l)$ is the right-invariant (left-invariant) Maurer-Cartan form.
    The group structure in $\sT G$ is given by the formula
        $$(X,g)\cdot (Y,h) = (X + \text{Ad}(g)Y, gh)
                                                                $$
     in the right trivialization and
        $$(g,X)\cdot (h,Y) = (gh, X + \text{Ad}(h^{-1})Y)
                                                                $$
in the left trivialization. The neutral element $e^{T}$ is represented by $(0,e)$
in the right trivialization and by $(e,0)$ in the left trivialization.

The Lie algebra $\sT \fr g$ of $\sT G$ is isomorphic as a linear space to $\fr
g\times \fr g$. This isomorphism is implemented by the zero section $j_G \colon G
\rightarrow \sT G$ and the obvious embedding  $j_{\fr g} \colon \fr g\rightarrow
\sT_e G \subset \sT G$ and is given by the following formula
        $$ \fr g\times \fr g \ni (X,Y) \rightarrow \sT_ej_G(X) + \sT_0 j_{\fr g}(Y)
\in \sT_{e^T}\sT G .
                                                    \tag 7.3        $$
 From now on we shall denote  $\sT_ej_G(X)$ by $\dot X$ and $\sT_0 j_{\fr g}(Y)$
by $\widehat{Y}$. We have the following commutation rules:
        $$\split [\dot X, \dot Y]_{\ssT \fr g} &= [X,Y]_{\fr g}^{\dsize \dot{}},\\
        [\dot X, \widehat{Y}]_{\ssT \fr g} = [\widehat{X},\dot Y]_{\ssT \fr g} &=
\widehat{[X,Y]}_{ \fr g}, \\
        [\widehat{X},\widehat{Y}]_{\ssT \fr g} &=0.
                                                \endsplit           $$
It follows that the Lie algebra $\sT \fr g$ is a semidirect product $\fr g
\ppr \fr g$ with respect to the adjoint representation in $\fr g$ of the Lie
algebra $\fr g$, i.~e.,
        $$ [(X,Y),(X',Y')]_{\fr g \ppr \fr g} = ([X,X']_{\fr g}, [X,Y']_{\fr g}
+[Y,X']_{\fr g}).
                                                                $$

    \proclaim{Theorem 7.2}{}
Let be $X\in \fr g$ and let $X^l_G$ be the corresponding left invariant vector field
on $G$. Then
        $$ \widehat{X}^l_{\ssT G} = (X^l_G)^v \ \ \text{and }\ \ \dot X^l_{\ssT G}
= (X^l_G)^c,
                                                                $$
i.~e., the corresponding to  $\widehat{X}$   and $\dot X$ left invariant vector
fields on $\sT G$ are the vertical and complete lifts of $X^l_G$ respectively. An
analogous statement is true for right invariant vector fields.
    \endproclaim
    \proof
The group multiplication $\sT m$ in $\sT G$ is the tangent lift of the group
multiplication $m$ in $G$. It follows that
        $$ \sT m(u_g,v_h) = \sT L_g(v_h) + \sT R_h (u_g),
                                                                $$
where $L_g$ and $R_h$ are left and right translations by $g$ and $h$ respectively.
The left translation $L^{\ssT}_{u_g}$ by $u_g$  in $\sT G$ is therefore given by
        $$L^{\ssT}_{u_g}(v_h) =  \sT L_g(v_h) + \sT R_h (u_g).
                                                                $$
     It is easy to verify that  $(X^l_G)^c(e^T) = \dot X$ and $(X^l_G)^v(e^{T})$.
Since the mapping  (7.3) is a linear  isomorphism, it is enough to show that vector
fields $(X^l_G)^c$ and $(X^l_G)^v$ are left-invariant on $\sT G$. In other words, we
have to show that flows they generate commute with left translations.

The flow $\psi^t$ of the vertical lift is given by
        $$\psi^t(v_h) = v_h + t X^l_G(h),
                                                                $$
      we have then
        $$ \split
L^{\ssT}_{u_g}\circ \psi^t(v_h) &= \sT L_g(v_h +tX_G^l(h)) + \sT R_h(u_g)\\
        &= \sT L_g(v_h) +\sT R_h(u_g) + t\sT L_g(X^l_G(h))\\
        &= \sT L_g(v_h) +\sT R_h(u_g) + t X^l_G(gh) = \psi^t \circ
L^{\ssT}_{u_g}(v_h).
                                        \endsplit                   $$
We made use of the left-invariance of $X^l_G$.

The flow  of  the complete lift $(X^l_G)^c$ is the tangent lift of the flow $\zf^t$
of $X^l_G$. We have
        $$\split L^{\ssT}_{u_g}\circ \sT\zf^t(v_h) &= \sT L_g \circ \sT\zf^t(v_h) +
\sT R_{\zf^t(h)}(u_g) \\
    &= \sT( L_g \circ \zf^t)(v_h) + \sT R_{\zf^t(h)}(u_g)
                                                \endsplit           $$
and, since $\zf^t$ is the flow of a left-invariant vector field, i.~e.,
        $$ L_g\circ \zf^t = \zf^t\circ L_g
                                                                $$
and
        $$\zf^t\circ R_h (g)=\zf^t(gh)= L_g(\zf^t(h)) = R_{\zf^t(h)}(g),
                                                                $$
 we get
        $$ \sT( L_g \circ \zf^t)(v_h) + \sT R_{\zf^t(h)}(u_g) = \sT( \zf^t \circ
L_g)(v_h) + \sT (\zf^t(h) \circ R_h) (u_g) = \sT\zf^t \circ L^{\ssT}_{u_g}(v_h).
                                                                $$

    \endproof

\subhead 8. Tangent Poisson-Lie algebras of Poisson-Lie groups  \endsubhead

    In this section we show that the tangent to the Poisson-Lie
algebra of a Poisson-Lie group $(G,\zL)$ is the Poisson-Lie
algebra of the Poisson-Lie group $(\sT G,\dt\zL)$. A special case
of of Poisson structures defined by $r$-matrices is discussed.

    The Lie algebra of a Poisson-Lie group inherits a Poisson structure. We recall here
a standard construction. More natural and more geometric one will be given in Section
10 (Proposition 10.4).

Let $(G,m,\zL)$ be a Poisson-Lie group. A Poisson structure  $\zL$ on a Lie group
can be regarded, using the right trivialization of $\sT G$, as a mapping  $\bar \zL
\colon G\rightarrow \fr g\wedge \fr g $. The Poisson structure $\zL$ is multiplicative
if and only if $\bar \zL$ is a
$1$-cocycle of $G$ with respect to the adjoint representation of $G$ in $\fr g\wedge
\fr g$.
 In particular, we have $\bar \zL (e) = 0$. The tangent mapping
        $$ \zl =\sT_e \bar \zL \colon  \sT_eG =\fr g \rightarrow \sT_0(\fr g
\wedge \fr g) = \fr g \wedge \fr g,
                                                                $$
being a 1-cocycle of $\fr g$,
     defines a cobracket (Poisson bracket on $\fr g$). The pair $(\fr g, \zl )$ is
called the {\it tangent Poisson-Lie algebra} of the Poisson-Lie group $(G,m, \zL)$.

Let $(X_1,\dots , X_n)$ be a basis of the Lie algebra $\fr g$. We can write
        $$\zL = \sum \zl^{ij} X^r_i\wedge X^r_j,
                                                                $$
where $\zl^{ij} $ are smooth functions on $G$ and $X^r_i$ are the
corresponding right invariant vector fields. We have then
        $$ \bar \zL(g) = \sum \zl^{ij}(g) X_i\wedge X_j
                                                                $$
and
        $$ \zl (X_k) =\langle X_k, \xd\zl^{ij}(e)\rangle X_i\wedge X_j
=\frac{1}{2}c^{ij}_kX_i\wedge X_j,
                                                                $$
where $c^{ij}_k$ are structure constants of the Lie algebra $\fr g$. The cobracket
$\zl \colon \fr g \rightarrow \fr g\wedge \fr g$ may be regarded as a bivector field
on $\fr g$ which defines a linear Poisson structure on $\fr g$.

    \proclaim{Definition}{} A {\rm Poisson-Lie algebra} $(\fr g,\zd )$ is a Lie
algebra $\fr g$ and a 1-cocycle (cobracket) $\zd \colon \fr g\rightarrow
\fr g \wedge \fr g$ (with respect to the adjoint representation of $\fr g$ in $\fr g
\wedge \fr g$)  such that it defines a Poisson structure on $\fr g$ or, equivalently,
that the dual mapping $ \zd^{\textstyle *} \colon \fr g^{\textstyle *} \wedge \fr
g^{\textstyle *} $ is a Lie bracket on $\fr g^{\textstyle *} $.
    \endproclaim The tangent Poisson-Lie algebra is an example of an abstract
Poisson-Lie algebra.  For  simple connected Lie groups there is one-to-one
correspondence between Poisson-Lie structures on $G$ and Poisson-Lie algebra
structures on $\fr g$. The Poisson-Lie algebra $(\fr g,\zl)$ can be seen as a
quadruple $(\fr g,\fr g^{\textstyle *}, [\cdot ,\cdot ],[\cdot ,\cdot ]^{\textstyle
*}  )$ and, for this reason, it is called sometimes a {\it Lie bialgebra}.

    \proclaim{Theorem 8.1}{}
Let  $\zl \colon \fr g\rightarrow \fr g\wedge \fr g $
 be the cobracket of the tangent Poisson-Lie algebra of a Poisson-Lie group $(G,m,\zL)$.
Then
        $$\dt \zl\colon  \sT\fr g \rightarrow \sT\fr g\wedge \sT\fr g
                                                                $$
is the cobracket of the tangent Poisson-Lie algebra of the tangent Poisson-Lie group
\newline $(\sT G,\sT m, \dt \zL)$. The dual mapping
        $$ (\dt\zl)^{\textstyle *} \colon (\sT \fr g)^{\textstyle *} \wedge (\sT
\fr g)^{\textstyle *} \rightarrow  (\sT \fr g)^{\textstyle *}
                                                                $$
is the Lie bracket on the tangent Lie algebra of $(\fr g^{\textstyle *}
,\zl^{\textstyle *} )$.
    \endproclaim
    \proof
    For any $X \in \zF^r(\ezt_M) $ and for any smooth function $f\in \zF^0(\ezp_M)$
we have
        $$ \dt(fX)= \ezt_M^{\textstyle *} (f)\dt X + \dt f\cdot \vt(X).
                                                                $$
Moreover, $Y^v(\ezt_M^{\textstyle *} f) =0$, $Y^v(\dt f) = Y^c(\ezt^{\textstyle *}
_M f) = \ezt^{\textstyle *} _M(Y(f))$ and $Y^c(\dt f) = \dt(Y(f))$ for any vector
field $Y\in \zF^1(\ezt_M)$.

Now, let  $\zL = \zL^{ij}X^r_i\wedge X^r_j$  and $\zl (X_k) = \frac{1}{2}
c^{ij}_k X_i\wedge X_j$  for a basis $(X_1,\dots,X_n)$ in $\fr g$. We have
        $$\split  \dt \zL  &= 2\ezt ^{\textstyle *} _G(\zL^{ij})(X^r_i
)^c \wedge (X_j^r)^v  + \dt(\zL^{ij})(X^r_i)^v \wedge (X^r_j)^v \\
        &= 2\ezt ^{\textstyle *} _G(\zL^{ij})(\dot X_i
)^r \wedge (\widehat{X}_j)^r  + \dt(\zL^{ij})(\widehat{X}_i)^r \wedge
(\widehat{X}_j)^r  .
                                                \endsplit       $$
     The cobracket $\zd$ on the Lie algebra of $\sT G$ is given by the formulae
        $$ \split  \zd (\widehat{X}_k) &= (X^r_k)^v(2\ezt_G^{\textstyle *}
(\zL^{ij})) (e^T)\dot X_i\wedge \widehat{X}_j + (X^r_k)^v (\dt(\zL^{ij}))(e^T)
\widehat{X_i}\wedge \widehat{X_j} \\
        &= \ezt_G^{\textstyle *} (X^r_k(\zL^{ij}))(e^T)\widehat{X}_i\wedge
\widehat{X}_j = X^r_k(\zL^{ij})(e)\widehat{X}_i\wedge \widehat{X}_j \\
        &= \langle X_k,\xd \zL^{ij}\rangle (e)\widehat{X}_i\wedge \widehat{X}_j =
\frac{1}{2}c^{ij}_k \widehat{X}_i\wedge \widehat{X}_j
                                            \endsplit               $$
and
        $$ \split \zd (\dot X_k ) &= (X^r_k)^c(2\ezt_G^{\textstyle *}
(\zL^{ij})) (e^T)\dot X_i\wedge \widehat{X}_j + (X^r_k)^c (\dt(\zL^{ij}))(e^T)
\widehat{X_i}\wedge \widehat{X_j} \\
    &= 2X^r_k(\zL^{ij})(e) \dot X_i\wedge \widehat{X}_j + \dt(X^r_j(\zL^{ij}))(e^T)
\widehat{X}_i\wedge \widehat{X}_j \\
    &= c^{ij}_k \dot X_i\wedge \widehat{X}_j
                                            \endsplit               $$
($\dt(f)$ is zero on the zero section).
It follows that $\zd = \dt \zl$.

    \endproof

Let $(G,m, \zL)$ be a Poisson-Lie group and let $(\fr g,[,],\fr p)$ be its Lie
bialgebra. Let us suppose that $\fr p$ is a coboundary (e. g.  $G$ is semisimple),
i.~e., that
        $$ \fr p (X) = [X,\r] = \r^{ij} ([X,X_i]\wedge X_j + X_i\wedge [X,X_j])
                                                                $$
for some $ \r= \r^{ij}X_i\wedge X_j \in \fr g\wedge \fr g$. It is known \cite{Dr}
that in this case the Poisson structure $\zL$ on $G$ can be written in the form
        $$ \zL = \r^l - \r^r,
                                                                $$
where $\r^l$ and $\r^r$ are the left- and right-invariant 2-vector fields
corresponding to $\r$. Since $\zL$ is a Poisson structure, $\r$ must satisfy the
generalized classical Yang-Baxter equation
        $$\text{ad}_X[\r,\r]=0
                                                                $$
for every $X\in \fr g$. The bracket  $[\r,\r ]$ is the algebraic Schouten bracket.
      An element of $ \fr g\wedge \fr g$ which satisfies this equation is called a
{\it generalized r-matrix} and the corresponding Poisson structure $\zL $ is called
{\it quasitriangular}.

    \proclaim{Theorem 8.2}{} Let $\zL = \r^l -\r^r$ be a quasitriangular Poisson-Lie
structure on a Lie group $G$ with the r-matrix $\r= \r^{ij}X_i\wedge X_j$, \ $\r^{ij}=
-\r^{ji}$. Then $\dt\zL$ is a quasitriangular Poisson-Lie structure on $\sT G$ with the
r-matrix $\dt \r = 2\r^{ij}\dot X_i\wedge \widehat{X}_j$.
    \endproclaim
    \proof
 Since $\zL = \r^{ij}(X^l_i\wedge X^l_j  - X^r_i\wedge X^r_j)$, we have
        $$\split \dt\zL &= 2\r^{ij}((X^l_i)^c\wedge (X^l_j)^v - (X^r_i)^c\wedge
(X^r_j)^v) \\
        &= 2\r^{ij}((\dot X_i)^l\wedge (\widehat{X}_j)^l -(\dot X_i)^r\wedge
(\widehat{X}_j)^r) = (\dt \r)^l - (\dt \r)^r
                                    \endsplit                       $$
and it is easy to check that $\dt\r$ is really an $\r$-matrix.
    \endproof

\subhead 9. Tangent lifts of generalized foliations  \endsubhead

    Symplectic  foliations of Poisson manifolds are
important examples of generalized foliations. In this section we
define the tangent lift of a generalized foliation and discuss
its basic properties.
    \proclaim{Definition}{}
    A {\it generalized distribution } on a manifold $M$ is a subset $S\subset \sT M$ such that
$S(x) = S\cap \sT_x M$  is a linear subspace for each point $x\in M$. $S$ is said to
be {\it smooth} if it is generated by the family
$$\Cal X (S) = \{X\in \zF^1(\ezt_M) \colon \forall x\in M \ \  X(x)\in S(x)\}
                                                                $$
of smooth vector fields, i.~e., $S(x)$ is spanned by $\{ X(x)\colon X\in \Cal X(S)\}$.

     A smooth distribution is {\it completely integrable} if for every point $x\in
M$ there exists an integral submanifold of $S$, everywhere of maximal dimension, which
contains $x$.

    The maximal integral submanifolds of  a completely integrable distribution $S$ form
a  partition of $M$, called the {\it generalized foliation} of
$M$ defined by $S$.
    \endproclaim
Let us notice that for a completely integrable distribution $S$ the family $\Cal
X(S)$ is a Lie subalgebra of the Lie algebra of vector fields on $M$.
The following theorem is due to H.~Sussmann \cite{Sus}.

    \proclaim{Theorem}{(Sussmann)} Let $S$ be a smooth distribution on $M$ and let $\Cal
D\subset \Cal X(S)$ be a family  of vector fields such that $\Cal D$ spans $S$. The
following properties are equivalent
    \roster
    \item $S$ is completely integrable,
    \item $S$ is invariant with respect to flows $\exp (tX)$ of vector fields
$X\in \Cal D$,
    \item flows of vector fields from $\Cal X(S)$ preserve $S$.
    \endroster
    \endproclaim

    \proclaim{ Theorem 9.1}{}
Let $S$ be a  completely integrable generalized distribution. Then the distribution
$S^{\ssT}$ generated by the family $\{ X^v,X^c\colon  X\in\Cal X(S)\}$ of vector fields
on $\sT M$ is completely integrable
    \endproclaim
    \proof
We have $\exp(tX^v)(v) = v + tX(\ezt_M v)$ and $\exp(tX^c) = \sT \exp(tX)$. Due to
the  formulae
        $$ (\exp(tX^v))_{\textstyle *}Y^v = Y^v,\ \  (\exp(tX^v))_{\textstyle
*}Y^c = Y^c + t[X,Y]^v
                                                                $$
     and
        $$ (\exp(tX^c))_{\textstyle *} Y^v =((\exp(tX))_{\textstyle *} Y)^v, (\exp
(tX^c))_{\textstyle *} Y^c =   ((\exp(tX))_{\textstyle *} Y)^c,
                                                                $$
it follows that $S^{\ssT}$ is invariant with respect to flows  of $X^c$ and $X^v$.
 From  the theorem of Sussmann, $S^{\ssT}$ is integrable.
    \endproof

    \proclaim{Definition}{}
    {\it The tangent foliation } $\Cal F^{\ssT}$ of a generalized foliation $\Cal F$
defined by $S$ is the foliation defined by $S^{\ssT}$.
    \endproclaim

    \specialhead Example. \endspecialhead

    Let us consider the distribution $S$ on $\R$, generated by
vector fields vanishing at $0\in \R$. The corresponding foliation
is the following one:
        $$\Cal F = \{ \{0\}, \R_+ ,\R_-\}, \ \  \text{where} \
\  \R_\pm = \{x\in\R \colon \pm x>0\}. $$
    We identify $\sT\R$ with $\R^2$ (with coordinates $(x,y)$)
and we get that vertical and complete lifts of vector fields
from $\Cal X (S)$ are of the form
        $$f(x) \frac{\partial }{\partial y}   $$
    and
        $$f(x)\frac{\partial }{\partial x} +
f'(x)y\frac{\partial }{\partial y}, $$
    where $f\in C^\infty(\R) $ vanishes at $0$.

    These vector fields generate the distribution $S^{\ssT}$
with
        $S^{\ssT}(x,y)= \text{span}\,(\dfrac{\partial
}{\partial x},\dfrac{\partial }{\partial y}) $ if $x\neq 0$,
$S^{\ssT}(0,y)= \text{span}\,(\dfrac{\partial }{\partial y}) $
if $y\neq 0$, and $S^{\ssT}(0,0) = \{ 0\}$.
Hence, the corresponding tangent foliation $\Cal F^{\ssT}$
consists of two half-plains $P_\pm =\{(x,y)\colon \pm x >0\}$,
two half-lines $L_\pm = \{(0,y)\colon \pm y>0\}$ and the point
$(0,0)$ as a 0-dimensional leaf.

    \proclaim{Proposition 9.2}{}
    If a 1-form $\zm$ annihilates a completely integrable distribution $S$ then
$\vt(\zm)$ and $\dt(\zm)$ annihilate $S^{\ssT}$.
    \endproclaim
    \proof
    Let  $X\in \Cal X(S)$. We have from Proposition 2.6 that
        $$ \langle \vt(\zm), X^v\rangle  =0, \ \ \langle \dt(\zm),X^c\rangle =
\dt\langle \zm,X\rangle =0
                                                                $$
and
        $$ \langle \vt(\zm), X^c\rangle = \langle \dt(\zm), X^v\rangle = \vt
(\langle \zm,X\rangle ) =0.
                                                                $$
    \endproof

    \proclaim{Proposition 9.3}{}
If a submanifold $N\subset M$ is a union of leaves of the foliation $\Cal F$ ($N$ is
$\Cal F$-foliated) then $\sT_NM = \ezt^{-1}_M(N)$  is $\Cal F^{\ssT}$-foliated.
    \endproclaim
    \proof It is enough to prove Proposition in the case of $N$ being a single
leaf.  Let $F$ be a leaf of $\Cal F^{\ssT}$. Since
$\sT\ezt_M(X^v)=0$ and $\sT\ezt_M(X^c)= X$ for any vector field $X$ on $M$, the tangent
projection  $\ezt_M(F)$  of $F$ is contained in a leaf of $\Cal F$. It follows that
$\sT_NM$ is $\Cal F^{\ssT}$-foliated.
    \endproof

    \proclaim{Proposition 9.4}{}
If $F$ is a leaf of $\Cal F$ then $\sT F$ is a leaf of $\Cal F^{\ssT}$.
    \endproclaim
    \proof
It is obvious that $\sT F\subset S^{\ssT} $ where $S$ is the generalized
distribution related to $\Cal F$. We have to show that $\sT F$ is maximal.
Since $F$ is maximal, $\sT _xF = S(x)$ and $S(F)$ is spanned by vector fields tangent
to $F$. On the other hand it is easy to see that if $X$ is a vector field on $M$,
tangent to $F$ on $F$, then $\dt X$ and $\vt X$ are tangent to $\sT F$ on $\sT F$.
Since $S^{\ssT}$ is generated by the family $\{ X^v, X^c \colon  X\in\Cal X(S)\}$,
$S^{\ssT}(\sT F) = \sT\sT F $, i.~e., $\sT F$ is an integral submanifold of $S^{\ssT}$
which is clearly  maximal.
    \endproof

\subhead 10. Symplectic foliations  of Poisson manifolds   \endsubhead

Let $(M,\zL)$ be a Poisson manifold. The characteristic distribution $S$ of $\zL$ is
generated by hamiltonian vector fields.  It is well known that $S$ is completely
integrable and that $\zL$ defines symplectic structures on leaves of $S$.
    \proclaim{ Proposition 10.1}{}
$S^{\ssT}$ is the characteristic distribution of $\dt\zL$.
    \endproclaim
    \proof
Since  the vertical and tangent lifts of 1-forms on $M$ generate the module of
1-forms on $\sT M$, it is enough to notice that (2.9) implies
        $$\langle \vt(\zm),\dt\zL\rangle = (\xi_{\zm} \zL)^v  \ \  \text{and}\ \
\langle \dt \zm, \dt\zL \rangle =(\xi_\zm \zL)^c
                                                                $$
and that, consequently, the characteristic distribution of $\dt\zL$ is generated by
the complete and vertical lifts of hamiltonian vector fields on $(M,\zL)$, i.~e., of
vector field from $\Cal X(S)$.
    \endproof

The following theorem by Weinstein \cite{We} describes the local structure of  a Poisson
manifold.

    \proclaim{Theorem}{(Weinstein)}
Let $(M,\zL)$ be a Poisson manifold of rank $2k$ at $x_0\in M$. Then there is an open
neighbourhood $U$ of $x_0$ such that $(U,\zl|_U)$ is isomorphic to a product
$(N\times V, \zL_N\times \zL_V)$ of Poisson manifolds where $(N,\zL_N)$ is a
symplectic manifold of dimension $2k$ and the rank of $(V,\zL_V) $ is zero at $z_0$,
$x_0 = (y_0,z_0)$.
    \endproclaim

The theorem of Weinstein shows that while analyzing only local properties of Poisson
manifolds it is enough to consider two cases:
    \roster
    \item $\zL$ is regular,
    \item $\zL$ vanishes at a point.
    \endroster

    \proclaim{Theorem 10.2}{}
Let $(M,\zL)$ be a Poisson manifold.
    \roster
    \item If $\zL$ is a regular Poisson structure of rank $2k$ then $\dt\zL$
is regular of rank $4k$.
    \item If $\zL$ is of rank 0 at $x_0\in M$ then $\sT_{x_0} M$ is a Poisson
submanifold of $(\sT M,\dt\zL)$ and $ \dt\zL$ defines on $\sT_{x_0}M$  a linear
Poisson structure (Kostant-Kirillov-Souriau structure). It induces then a Lie algebra
structure on $\sT_{x_0}^{\textstyle *}M $.
    \endroster
    \endproclaim
    \proof
    \roster
    \item It follows from the theorem of Weinstein that we can choose local coordinates
on $M$ such that
        $$ \zL = \sum_{j,i=1}^r \zL^{ij}\frac{\partial }{\partial x^{i}}\wedge
\frac{\partial }{\partial x^{j}},\ \ \ \zL^{ij} =- \zL^{ji}, \ \ \det (\zL^{ij})\neq
0,
                                                                $$
where $\zL^{ij}$ are constant.
Hence, in the adopted coordinate system,
        $$ \dt\zL = 2\sum_{j,i=1}^r \zL^{ij}\frac{\partial }{\partial \dot x^{i}}
\wedge \frac{\partial }{\partial x^{j}}.
                                                                $$
     \item We have locally
        $$ \zL (x) = \sum_{j,i} \zL^{ij}(x)\frac{\partial }{\partial x^{i}}\wedge
\frac{\partial }{\partial x^{j}},\ \ \ \zL^{ij} = -\zL^{ji}, \ \ \zL^{ij}(x_0)= 0
                                                                $$
and
    $$\dt\zL (x_0, \dot x) = \sum_{j,i,k} \frac{\partial \zL^{ij} }{\partial
x^{k}}(x_0)\dot x^k \frac{\partial }{\partial \dot x^{i}}\wedge \frac{\partial
}{\partial\dot x^{j}}.
                                                                $$
Hence $\sT_{x_0} M$ is a Poisson submanifold with the Poisson bracket
        $$\{ \dot x^i, \dot x^j \} = 2\sum_k \frac{\partial \zL^{ij} }{\partial
x^{k}} (x_0)\dot x^k.
                                                                $$
    \endroster
    \endproof
    \proclaim{Corollary}{}
If $\zL$ is of rank 0 at $x_0$ then the Poisson structure $\dt\zL$ on $\sT_{x_0}$ is
given by the formula
        $$ \dt \zL (v) = \vt(\ll_{\widetilde{v}}\zL)(v) , \ \ v\in\sT_{x_0}M,
                                                                $$
where $\widetilde{v}\in \zF^1(\ezt_M)$ is such that $\widetilde{v}(x_0)=v$.
    \endproclaim
    \proof
Since $\zL(x_0) = 0$ it follows that $\dt\zL(v)$ is vertical and,
consequently,
$$\vt(\widetilde{v}^{\textstyle *} \dt\zL)(v) = \dt\zL(v) $$ (see Section 3). It
follows from Theorem 3.2 that
        $$ \dt\zL(v) = \vt(\ll_{\widetilde{v}}\zL)(v).
                                                                $$
    \endproof

The following two propositions complete our discussion on the structure of the
tangent Poisson manifold.

    \proclaim{Proposition 10.3}{}
Let $f$ be a local Casimir of $\zL$, i.~e., $\zL(\xd f,\cdot ) =0$. Then $\vt(f)$ and
$\dt (f)$ are local Casimirs of $\dt\zL$.

Moreover, if $\zL$ is regular at $x\in M$ with symplectic leaves determined by
Casimirs $(f_1, \dots, f_n)$, then $\zL$ is regular at $v\in \sT_xM$ with  symplectic
leaves determined by Casimirs $(\ezt_M^{\textstyle *}( f_1),
\dots,\ezt_M^{\textstyle *}( f_n) ,\dt f_1, \dots, \dt f_n)$.
    \endproclaim
    \proof
It follows from (2.9) that if $f$ generates a zero-hamiltonian field then also
$\ezt_M^{\textstyle *} f$ and $\dt f$ generate the zero-hamiltonian field on $\sT M$,
i.~e., they are Casimirs.

Let $x\in M $ be a regular point of $\zL$  with the symplectic foliation in the
neighbourhood of $x$ defined by Casimirs $(f_1,\dots,f_n)$. We may assume that $\xd
f_1, \dots , \xd f_n$ are linearly independent at $x$. It follows that
$(\xd\ezt_M^{\textstyle *}( f_1), \dots,\xd \ezt_M^{\textstyle *}( f_n) ,\xd \dt
f_1, \dots, \xd \dt f_n)$ are linearly independent at $v\in \sT_x M$. Since the
rank of $\dt\zL$ is $2(\dim M -n)$ (Theorem 10.2) the theorem follows.
    \endproof

If $(G,\zL)$ is a Poisson-Lie group then $\zL(e) =0$ and identifying $\fr g $ with
$\sT_e G$ we get a Poisson structure on $\fr g$ induced by $\dt\zL$.
    \proclaim{Proposition 10.4}{}
The Poisson structure on $\fr g$ induced by $\dt\zL$ is equal to the Poisson
structure defined by the cobracket of the tangent Poisson-Lie  algebra.
    \endproclaim
    \proof
Let $(x^i)$ be a coordinate system on $G$ centered at $e$, i.~.e, $x^i(e)=0$.
It defines a basis $(X_i)$ of the Lie algebra $\fr g$, $X_i = \dfrac{\partial
}{\partial x^i}(e)$.  There are functions $a^k_i$ on $G$ such that
        $$\frac{\partial }{\partial x_{i}} - X_i^r = a_i^k
X_i^r, $$
    where $a_i^k(0)=0$.
We have then
        $$\zL = \zL^{ij} \frac{\partial }{\partial x^{i}}\wedge
\frac{\partial }{\partial x^{j}} = \zl^{ij} X^r_i\wedge X^r_j
                                                                $$

and
        $$ \zL^{ij}(0) = \zl^{ij}(0).
                                                                $$
The induced  cobracket $\zl$ is given by
        $$\zl(X_k) = \dfrac{\partial \zl^{ij} }{\partial x_{k}}(0)X_i\wedge X_j.
                                                                $$
    On the other hand,
        $$\dt|_{\ssT_eG} = \frac{\partial \zL^{ij} }{\partial x_{k}}(0)\frac{\partial
}{\partial x_{i}}\wedge \frac{\partial }{\partial x_{j}} = \frac{\partial \zL^{ij}
}{\partial x_{k}}(0)X_i\wedge X_j.
                                                                $$
  It follows that
        $$ \zl^{ij}= \zL^{ij} + \zL^{il}a^j_l - \zL^{jl}a^i_l
+\frac{1}{2}(a^i_ka^j_l - a^j_k a^i_l)\zL^{kl}.
                                                                $$
Hence
        $$ \dfrac{\partial \zl^{ij} }{\partial x_{k}}(0) = \dfrac{\partial
\zL^{ij} }{\partial x_{k}}(0).
                                                                $$
    \endproof

\subhead 11. Examples \endsubhead
\specialhead Example 1 \endspecialhead
    On $\fr{su}(2)^{\textstyle *} \simeq \R^3$ consider the linear Poisson structure
        $$ \zL = z\frac{\partial }{\partial x}\wedge \frac{\partial }{\partial y} +
x\frac{\partial }{\partial y}\wedge \frac{\partial }{\partial z} +
y\frac{\partial }{\partial z}\wedge \frac{\partial }{\partial x}.
                                                                $$
The symplectic foliation of $\R^3$ is the union of 2-dimensional spheres $x^2 + y^2
+z^2 =r>0$ and the origin $(0,0,0)$ as a 0-dimensional leaf. It is regular outside
the origin and is defined by the Casimir $f(x,y,z) = x^2 +y^2 +z^2$.

The tangent Poisson structure is given by the formula
        $$ \multline \dt\zL = z(\frac{\partial }{\partial \dot x}\wedge
\frac{\partial }{\partial y} + \frac{\partial }{\partial x}\wedge \frac{\partial
}{\partial \dot y} ) + x(\frac{\partial }{\partial \dot y}\wedge \frac{\partial
}{\partial z} + \frac{\partial }{\partial y}\wedge \frac{\partial }{\partial \dot
z}) + y(\frac{\partial }{\partial \dot z}\wedge \frac{\partial }{\partial x} +
\frac{\partial }{\partial z}\wedge \frac{\partial }{\partial \dot x})\\
    + \dot z \frac{\partial }{\partial \dot x}\wedge \frac{\partial }{\partial \dot
 y} + \dot x \frac{\partial }{\partial \dot y}\wedge \frac{\partial }{\partial \dot
 z} + \dot y \frac{\partial }{\partial \dot z}\wedge \frac{\partial }{\partial \dot
x} .
                                        \endmultline                $$
The symplectic foliation of $\sT \R^3$ is regular outside $\sT_0 \R^3$ and it is
determined by two Casimirs
        $$ f_0(x,y,z,\dot x,\dot y,\dot z) = x^2 + y^2 + z^2 , \ \ f_1(x,y,z,\dot
x,\dot y,\dot z) = x\dot x + y\dot y + z\dot z .
                                                                $$
The tangent space $\sT_0 \R^3 \simeq \R^3$ has a linear Poisson structure
        $$ \dot z\frac{\partial }{\partial \dot x}\wedge \frac{\partial }{\partial
\dot y} + \dot x\frac{\partial }{\partial \dot y}\wedge \frac{\partial }{\partial
\dot z} + \dot y\frac{\partial }{\partial \dot z}\wedge \frac{\partial }{\partial
\dot x}.
                                                                $$
which is equal to $\zL$.

\specialhead Example 2 \endspecialhead
In this example $M=\R^4$ and $\zL$ is the following quadratic Poisson structure:
        $$\multline \zL = x_1x_3\frac{\partial }{\partial x_{3}}\wedge
\frac{\partial }{\partial x_{2}} +x_1x_4\frac{\partial }{\partial x_{4}}\wedge
\frac{\partial }{\partial x_{2}} + x_2x_3\frac{\partial }{\partial x_{1}}\wedge
\frac{\partial }{\partial x_{3}} \\ + x_2x_4\frac{\partial }{\partial x_{1}}\wedge
\frac{\partial }{\partial x_{4}} + (x_3^2 + x_4^2) \frac{\partial }{\partial
x_{2}}\wedge \frac{\partial }{\partial x_{1}}.
                                \endmultline                        $$
This structure is degenerated at points of the linear subspace $x_3=x_4 =0$ and regular
outside it. The symplectic foliation is the intersection of the ,,book'' foliation
generated by vector fields $x_3\dfrac{\partial }{\partial x_{3}} +x_4\dfrac{\partial
}{\partial x_{4}} $, $\dfrac{\partial }{\partial x_{1}}$ and $\dfrac{\partial
}{\partial x_{2}}$, consisting of 3-dimensional half-spaces sewed up along the
2-dimensional edge $x_3=x_4=0$ of 0-dimensional leaves, and the spherical foliation.

The unit sphere $S^3= \{x\colon \sum x_i^2 =1 \}$ we identify with the $SU(2)$ group by
        $$ (x_1,x_2,x_3,x_4) \longmapsto \pmatrix x_1 + i x_2 & -x_3 + i x_4 \\
            x_3 + i x_4 & x_1 - i x_2 \endpmatrix .
                                                                $$
With this identification, the tensor $\zL$ restricted to $S^3$ defines a Poisson-Lie structure
on $SU(2)$. This structure is quasitriangular with the generalized r-matrix
        $$ \r = \frac{\partial }{\partial x_{3}}\wedge \frac{\partial }{\partial
x_{4}} \in \fr{su}(2)\wedge \fr{su}(2).
                                                                $$
      Here we used identifications
        $$ \fr{su}(2) = \sT_eSU(2) \simeq \sT_{(1,0,0,0)} S^3,
                                                                $$
so in the matrix form,
        $$ \r = X\wedge Y \ \ \text{with} \ \  X = \pmatrix 0 & -1 \\ 1 & 0
\endpmatrix,\ \ Y = \pmatrix 0 & i \\ i & 0 \endpmatrix \in \fr{su}(2).
                                                                $$
     It is worthy  noticing that this Poisson-Lie structure is related to the
quantum $SU(2)$ group of Woronowicz (cf. [Gr]).

Singular points on $S^3$ form a 1-dimensional circle    $x_1^2 + x_2^2 =1,\ x_3 =x_4
=0$ which corresponds to the Cartan subgroup (maximal torus) of $SU(2)$.
At a singular point $x= (\cos \zf,\sin \zf,0,0)$, the tangent space $\sT_xS^3$ carries
a linear Poisson structure induced by $\dt \zL$
        $$\dt\zL|_{\ssT_xS^3} = \cos (\zf)(\dot x_4 \frac{\partial }{\partial
\dot x_{4}}\wedge \frac{\partial }{\partial \dot x_{2}} + \dot x_3 \frac{\partial }{\partial
\dot x_{1}} \wedge \frac{\partial }{\partial \dot x_{2}}) + \sin (\zf)(\dot x_3
\frac{\partial }{\partial \dot x_{1}} \wedge \frac{\partial }{\partial \dot x_{3}} +
\dot x_4 \frac{\partial }{\partial \dot x_{1}} \wedge \frac{\partial }{\partial \dot x_{4}})
                                                                $$
associated with the Lie bracket
        $$ \split \{ \dot x_1, \dot x_3\} &= \sin(\zf) \dot x_3 , \qquad \\ \{ \dot
x_1, \dot x_4\} &= \sin(\zf) \dot x_4 ,\\
    \{ \dot x_2, \dot x_3\} &= -\cos(\zf) \dot x_3, \qquad  \\ \{ \dot x_2, \dot
x_4\} &= -\cos\zf) \dot x_3  .
                                    \endsplit                       $$
In particular, on $\sT_eSU(2) = \sT_{(1,0,0,0)}S^3$ we have
        $$ \dot x_4 \frac{\partial }{\partial \dot x_4} \wedge \frac{\partial
}{\partial \dot x_2} + \dot x_3 \frac{\partial }{\partial \dot x_3}\wedge
\frac{\partial }{\partial \dot x_2}.
                                                                $$

It follows that the cobracket $\zl$ of the tangent bialgebra of the Poisson-Lie group
$SU(2)$ is given by
        $$ \zl(\frac{\partial }{\partial \dot x_4}) = \frac{\partial }{\partial
\dot x_4}\wedge  \frac{\partial }{\partial \dot x_2}, \ \  \zl(\frac{\partial
}{\partial \dot x_3}) = \frac{\partial }{\partial \dot x_3}\wedge \frac{\partial
}{\partial \dot x_2}, \ \ \zl(\frac{\partial }{\partial \dot x_2}) = 0
                                                                $$
and the associated Lie bracket on $\fr{su}(2)^{\textstyle *} $ is given by
        $$\{\dot x_4, \dot x_2\} =\dot x_4,\ \  \{\dot x_3, \dot x_2\} = \dot x_3,
\ \ \{ \dot x_4, \dot x_3\} = 0.
                                                                $$

We recognize this structure as the structure of the Lie algebra $\fr{sb} (2,\C)$ of
$2\times 2$ traceless, upper triangular complex matrices with real diagonal elements.

The introduced Poisson-Lie structure on $SU(2)$ defines then the
group $SB(2,\C)$ as the dual group. It is not difficult to verify
that the corresponding double group can be identified with
$SL(2,\C)$ with $SU(2)$ and $SB(2,\C) $ canonically embedded as
subgroups.
\subhead 12. Conclusions \endsubhead

In this paper we introduced  lifts of multivector fields
and related objects (like generalized foliations) from a manifold
$M$ to its tangent bundle. These operations can be considered as
an extention of the tangent functor to these objects and
corresponding structures. We called them tangent lifts.

Among new results are those stating that the tangent lift
operator $\dt$ on multivector fields commutes with the Schouten
bracket (Theorem 2.4), that the symplectic foliation of the
tangent Poisson structure is the tangent foliation of the given
Poisson structure (Proposition 10.1), and that the tangent lift
of a Poisson-Lie structure is a Poisson-Lie structure (Theorem
7.1), etc. We proved also  theorems describing Poisson structures by
conditions for morphisms of the tangent and cotangent bundles
(Theorem~4.4 and Theorem~4.5).

Some of results refer to already known facts, but the
used methods give us new point of view, show better relations
between different objects and provide deeper understandig of
some well-known structures and facts as special cases of more
general situations (cf. Theorem~5.1 and Theorem~10.2).

We are convinced that these results  show the importance of the
concept of tangent structures in general and of the derivation
$\dt$ in particular.   The next step should be the analysis of
multitangent constructions,  important for classical field
theory (multiphase approach) and classical mechanics, of
extended objects. As we have seen, the conditions discussed in
Section 4 are, in fact, compatibility conditions of tangent and
canonical structures. This idea can be applied in more general
situations like in nonrelativistic, time dependent mechanics,
where the structure needed is more general than Poisson or
symplectic  one (\cite{Ur}).  Results of further studies on tangent
lifts together with applications to analytical mechanics will be
given in a separate publication.

\subhead Acknowledgments \endsubhead

This work is a contribution to a program of geometric formulation
of phyasical theories conducted jointly with W.~M.~Tulczyjew. The
authors wish to thank G.~Marmo for helpful suggestions concerning
the case of Poisson-Lie structures.

\Refs \widestnumber \key{MFLMR}
     \blb By \by G.~Byrnes \paper A complete set of Bianchi
     identities for tensor fields along the tangent bundle projection
     \jour J. Phys. A: Math. Gen. \vol 27 \yr 1994 \pages 6617--6632
     \endbb
     \blb Co1 \by T.~Courant \paper Tangent Dirac structures \jour J.
     Phys. A: Math. Gen. \vol 23 \yr 1990\pages 5153--5168 \endbb
    \blb Co2 \by T.~Courant \paper Tangent Lie Algebroids \jour J.
     Phys. A: Math. Gen. \vol 27(13) \yr 1994\pages 4527--4536    \endbb
    \blb Co3 \by T.~Courant \paper Tangent Poisson Structure \paperinfo unpublished\endbb
    \blb Dr \by V.~G.~Drinfeld\paper Hamiltonian structures on Lie groups, Lie
bialgebras and the geometric meaning of the classical Yang-Baxter equation \jour
Soviet Math. Dokl.\yr 1983\vol 27(1)\pages 68-71 \endbb
    \blb Gr \by J.~Grabowski \paper Quantum $SU(2) $ group  of
Woronowicz and Poisson structures \inbook Differential Geometry
and Its Applications \eds J.~Jany\v{s}ka and D.~Krupka \bookinfo
Proc. Conf. Brno (1989) \publ World Scientific \yr 1990 \pages
313-322 \endbb

    \blb KMS \by I.~Kol\'a\v{r}, P.~W.~Michor and
J.~Slovak \book Natural Operations in Differential Geometry
\publ Springer-Verlag  \yr 1993 \endbb

    \blb MCS1 \by
E.~Mart\'{\i}nez, J.~F.~Cari\~nena and W.~Sarlet  \paper
Derivations of differential forms along the tangent bundle
projection \jour Diff. Geom. Applic. \vol 2 \yr 1992 \pages
17--43  \endbb
    \blb MCS2 \by E.~Mart\'{\i}nez, J.~F.~Cari\~nena
and W.~Sarlet  \paper Derivations of differential forms along
the tangent bundle  projection II\jour Diff. Geom. Applic. \vol
3 \yr 1993 \pages 1--29  \endbb
    \blb MCS3 \by E.~Mart\'{\i}nez,
J.~F.~Cari\~nena and W.~Sarlet  \paper Geometric
characterization of separable second-order  differential
equations \jour Math. Proc. Camb. Phil. Soc. \vol 113  \yr 1993
\pages 205--224 \endbb
    \blb Mi \by P.~W.~Michor \paper Remarks
on the Schouten-Nijenhuis  bracket \jour Suppl. Rend. Circ. Mat.
Palermo, Serie  II\vol 16\yr 1987\pages 207--215\endbb
    \blb MFL \by G.~Morandi, C.~Ferrario, G.~ Lo Vecchio,
G.~Marmo and C.~Rubano\paper The inverse problem in the calculus
of variations  and the geometry of the tangent bundle \jour
Physics Reports \vol 188  \yr 1990 \pages 147--284\endbb
    \blb PiTu \by G.~Pidello  and W.~M.~Tulczyjew\pages 249--265
\paper Derivations of differential forms on jet bundles \yr1987
\vol 147\jour Ann.Mat.Pura Appl.\endbb
\blb Sus \by H.~J.~Sussmann \pages 171--188 \paper Orbits of families
of vector fields and integrability of distributions \yr1973
\jour Trans. Amer. Math. Soc. \vol 180 \endbb
    \blb SdA \by G.~S\'anchez de Alvarez \pages 399--409 \paper
Controllability of Poisson Control Systems with Symmetries \jour
Contemporary Mathematics \vol 97 \yr 1989 \endbb
    \blb Tu1 \by W.~M.~Tulczyjew \pages 101--114 \paper
Hamiltonian systems, Lagrangian systems, and the Legendre
transformation \jour Ann.~Inst.~H.~Poincar\'e \yr1977\vol27\endbb
    \blb Tu2 \by W.~M.~Tulczyjew  \paper A symplectic formulation of relativistic
particle dynamics \jour Acta Phys. Pol.  \vol {\bf B }8 \yr 1977 \pages
431--447\endbb
    \blb Ur \by P.~Urba\'nski \paper Affine Poisson structures in analytical
mechanics \inbook Quantization and Infinite Dimensional Systems
\eds J.--P.~Antoine and A.~Odzijewicz\yr 1994 \pages 123--129
\publ Plenum Press, New York and London  \endbb
    \blb We \by A.~Weinstein \pages 523--557 \paper The local structure of Poisson
manifolds \yr1983\jour J. Differential Geometry \vol 22\endbb
\endRefs

\enddocument